\documentclass[11pt]{article}
\usepackage{graphicx}
\usepackage{amsmath, bbm}\usepackage{mathrsfs,amsfonts, bbm}
\usepackage{amsfonts}
\usepackage{amsmath,amssymb,amsfonts}
\usepackage{amssymb}
\usepackage{color}
\usepackage{cases}
\definecolor{c}{rgb}{0.9,0.3,0.1}
\definecolor{b}{rgb}{0.1,0.3,0.9}
 
\allowdisplaybreaks

\newtheorem{theorem}{Theorem}[section]

\newtheorem{lemma}[theorem]{Lemma}
\newtheorem{prop}[theorem]{Proposition}
\newtheorem{definition}[theorem]{Definition}

\newtheorem{remark}[theorem]{Remark}

\newtheorem{hypothesis}[theorem]{Hypothesis}

\def\al{{\alpha}}\def\be{{\beta}}
\def\ep{{\varepsilon}}\def\ka{{\kappa}}
\def\la{{\lambda}}
\def\si{{\sigma}}\def\va{{\varepsilon}}

\def\<{\langle}\def\>{\rangle}

\def\({\left(}\def\){\right)}

\newfam\msbmfam\font\tenmsbm=msbm10\textfont
\msbmfam=\tenmsbm\font\sevenmsbm=msbm7
\scriptfont\msbmfam=\sevenmsbm

\def\EE{ \mathbbm E}\def\PP{\mathbbm P}
\def\RR{\mathbbm R}

\def\cF{{\cal F}}
\def\cM{{\cal M}}

\def\cU{{\cal U}}

\newcommand{\bP}{{\mathbbm{P}}}
 \newcommand{\E}{{\mathbbm{E}}}

\newcommand{\R}{{\mathbbm{R}}}

\def\R{\mathbbm{R}}
\def\1{\mathbbm{1}}

 \numberwithin{equation}{section}
 
\vspace{4ex}

\def\qed{\hfill $\Box$}

\def\sF{\mathcal F}

\def\sH{\mathcal H}

\def\wt{\widetilde}
\def\wh{\widehat}

\def\eps{\varepsilon}

\def\pf{\noindent {\bf Proof.} }

\begin{document}

\title{\huge
Stochastic maximum principle   for sub-diffusions  and its applications
\thanks{This is the arXiv version,   containing additional  details,  of the paper  with the same title that is going to appear in 
{\it SIAM J. Control Optim. }}
}
 \author{ Shuaiqi Zhang\footnote{Research partially  supported  by the Fundamental Research Funds for the Central Universities (Grant No. 2020ZDPYMS01) 
 and Nature Science Foundation of Jiangsu Province (Grant No. BK20221543)  and National Natural Science Foundation of China (Grant No. 12171086)   .}   \\ \small
 School of Mathematics, China  University of Mining and Technology,
\\ \small
Xuzhou, Jiangsu, 221116, China\\
\\
 {Zhen-Qing Chen\footnote{Corresponding author.  Research  partially supported by a Simons Foundation Grant. }
 }\\  \small Department of Mathematics,  University of Washington,
\\
 \small Seattle, WA 98195, USA\\
    }
\maketitle
\date

\begin{abstract}
In this paper, we study optimal stochastic  control  problems  for stochastic systems driven by non-Markov sub-diffusion $B_{L_t}$,
which have the mixed features of   deterministic  and   stochastic  controls.
 Here  $B_t$ is the  standard Brownian motion on $\RR$, and
  $L_t:= \inf\{r>0: S_r>t\}, \quad t\geq 0,$
  is the inverse of  a subordinator  $S_t$ with drift $\kappa >0$ that is independent of $B_t$.
 We obtain  stochastic maximum principles (SMP) for these systems using both  convex  and  spiking variational  methods,
 depending on whether the convex domain is convex or not.   
  To derive SMP,  we first establish a  martingale representation theorem for sub-diffusions $B_{L_t}$,
and  then  use it to derive the  existence and uniqueness result for the solutions of  
 backward stochastic differential equations (BSDEs) driven by sub-diffusions, which may be of independent interest. 
  We also derive sufficient SMPs. 
   Application to a linear quadratic system is   given to illustrate the main results of this paper.
   
\medskip

 \noindent {\bf Keywords}: Sub-diffusion, 
 stochastic maximum principle, BSDE driven by sub-diffusion, martingale representation theorem of sub-diffusion 
  
\medskip

\noindent{\bf AMS 2020 Subject Classification:} Primary 93E20, 60H10; 
 Secondary 49K45
\end{abstract}

\section{Introduction}\label{intro}

This paper is concerned with  optimal stochastic  controls for stochastic differential  equations (SDEs) driven by anomalous 
 sub-diffusions,
which are the time change of Brownian motions by inverse subordinators. 
  Sub-diffusions  are  random  processes that  describe the motions of particles that moves slower than Brownian motion. 
 They       have  been widely used to model  many natural systems
   \cite{KS,   MK}, 
 ranging
 from   transport processes in porous media 
 as well as in systems with memory  
  to avascular tumor growth.    Unlike Brownian motion, anomalous sub-diffusion by itself is not a Markov process.

The study of   sub-diffusion   is also an active research area in mathematics. 
 There are many work on sub-diffusions and the corresponding fractional time differential equations.
See, for example, 
 \cite{Meerschaert2004, Meerschaert2014}
  on   sub-diffusion processes as  scaling limits of continuous time random walks with heavy tail waiting times, 
 \cite{Baeumer2009, Chen2017, CKKW1, CKKW2}  on  the connections between  the sub-diffusions and the time fractional   equations  and on the estimates of the fundamental solutions,  and 
\cite{Marcin}  on some  sample path properties of sub-diffusions.
 In this paper we study   stochastic control problems for systems driven by anomalous sub-diffusions.  
  As one will see from below, it is quite natural to consider such stochastic systems as they can be used to model 
  bear markets in which activities are much slower.

We now describe the setting and the model in more details. 
Suppose that $S_t$ is a subordinator with drift $\kappa >0$ and  
L\'evy measure $\nu$; that is, $S_t=\kappa t+ S_t^0$, where $S_t^0$ is a driftless subordinator with
  L\'evy measure $\nu$.
 Let $ L_t:= \inf\{r>0: S_r >t\}$, $t\geq 0$, be the inverse of $S$. 
 The inverse subordinator $L_t$ is continuous in $t$  but stays constant during infinitely many time periods which are resulted from the infinitely many jumps by the subordinator $S_t$ when its L\'evy measure is non-trivial.  
   Clearly,   $0\leq L_{t+s} -L_t \leq s/\kappa$ for any for any $t, s>0$.
 So almost surely
\begin{equation}\label{e:1.1}
K_t:=\frac{dL_t}{dt}    \ \hbox{ exists} \quad \hbox{and} \quad 0\leq K_t \leq 1/\kappa   \    \hbox{ for a.e. }   t>0.
\end{equation} 
Let $B$ be a  Brownian motion independent of $S_t$. The sub-diffusion $B_{L_t}$ is a continuous martingale with
$\langle B_{L_t} \rangle = L_t$.
  In this paper,  we study stochastic optimal control for systems driven by the anomalous sub-diffusion $B_{L_t}$.

  One of our motivations comes from 
 stochastic control problems in finance such as optimal investment problems.  
Black-Scholes model is typically employed to describe the   stock price, i.e.,
 $
 dS_t=  \ \mu_t  S_t  dt +\sigma_t  S_t  d B_t$,  where $\mu_t$ is the stock return  rate, $\sigma_t$ is the volatility, and  $B_t$ is a standard Brownian motion.
 To model   the bear  and bull markets,   most papers use 
  a  regime switching Markov chain between   different values for  the return rate   function  $\mu_t$,
   that  is,   one uses 
 use a large return value for $\mu_t$ during  the bull market and  a smaller return value for $\mu_t$  during the bear market;  see, e.g.,  \cite{zhang2018}. However this kind of 
 Markov regime switching models   do not  capture the phenomenon that tradings are typically much less active
 in the bear market. We propose to    use 
sub-diffusion $B_{L_t}$ in place of Brownian motion $B_t$ to model the bear market:
  $\displaystyle          dS_t=\mu_t S_t dt +\sigma_t S_t d B_{L_t}$. 
  Note that  $B_{L_t}$ stays flat during the time periods when $L_t$ stays constant. When the L\'evy measure $\nu$ of  the subordinator $S_t$ is infinite,
   $S_t$ has infinitely many small jumps in any given time interval and, for any $\eps >0$,  the jumps of $S_t$ of size larger than $\eps$ occurs 
   according to a Poisson process with parameter $\nu (\eps, \infty)$.  In this case, during any fixed time intervals, 
   $L_t$ has infinitely many small time periods    but  only finite many time periods larger than a fixed length  during which it stays constants.
   Thus the sub-diffusion $B_{L_t}$  matches well the phenomena that the market constantly  has small corrections 
    but big corrections or bear market occurs only sporadically.

  Stochastic maximum principle (SMP in short), which gives necessary conditions for   a stochastic control to be optimal, is a
fundamental principle for  optimal  stochastic  controls.  
  In  Brownian models (that is, when the subordinate $S_t$ is just the deterministic process $\kappa t$ and 
	  thus  $B_{L_t}$ is a Brownian motion), 
 it has  been studied extensively since the pioneering work of   Kushner \cite{Kushner1965, Kushner1972}. 
  In \cite{Peng1990}, Peng   introduced the second-order term in the Taylor expansion of the variation 
and obtained the global maximum principle for the stochastic optimal control problem. 
Pardoux and Peng \cite{PardouxPeng1990} established a strong connections between  
  backward stochastic differential equations (BSDEs in abbreviation) 
 and stochastic control problems.
Since then, many researchers  have investigated 
 optimal control problems for various stochastic systems
(see, for example,   \cite{Tang1994,   Oksendal2005,  Li2012,  zhang2020}).
   When the control domain  is not convex, the   spiking variation method is  employed to  study  SMP,  while
when the control domain is convex, convex variation method is used    for    SMP.

    In this paper,  we study   and derive stochastic maximum principles for SDEs \eqref{e:state}
    driven by sub-diffusions 
   using   both the convex variational method and the spiking variational method
depending on whether the control domain is convex or not. 
   Since the subordinator $S_t$ contains the discontinuous component $S^0$,  many new phinomena and difficulties arise.
    It seems that this is the first time  stochastic maximum principles for systems driven 
    by sub-diffusions have been systematically investigated. 
   We further establish two sufficient stochastic maximum principles (that is, sufficient conditions for a stochastic control to be optimal)  for SDEs driven by sub-diffusion, one for the general case and the other for the convex control domain case.
   See Remark \ref{R:1.1} below for a recent related work by Nane and Ni \cite{NN} on sufficient SMP.  
   The novelties and main contributions of this paper  are as  follows.
   
   \smallskip
   
 \begin{enumerate}
\item[1)]  The L\'evy measure $\nu$ of the subordinator $S_t$ in this paper can be any L\'evy measure, finite or infinite.  
 SDEs driven by such anomalous sub-diffusions can be used to model variety of situations such as bull-and-bear markets.
Note that when $\nu$ is a non-zero finite measure,  $S_t$ is a compounded Poisson process with positive drift $\kappa$.
Thus ignoring the flat time duration, the arrivals of the flat time interval of $B_{L_t}$ (i.e., inactive time period) is Poisson 
and the length of the flat time durations are determined by the corresponding jump size of $S_t$ multiplied by $1/\kappa$. 
When the L\'evy measure $\nu$ is infinite, in any given time interval there are infinitely many tiny inactive time periods for $B_{L_t}$ 
but for any $\eps >0$, the number of inactive time periods larger than $\eps$ is roughly Poisson. 
When $\nu =0$, $B_{L_t}$ reduces to a Brownian motion. Thus the results of this paper extend the corresponding 
classical stochastic optimization results in a continuous and stable way: when $\nu=0$, they recover the corresponding classical results
in the Brownian setting. 

\item[2)] In the typical case that the L\'evy measure $\nu$ is non-zero, 
the control problem is   \emph{not   entirely    stochastic} as there are intervals on which the sub-diffusion $B_{L_t}$ is constant.
On the other hand, these flat time intervals are random not deterministic. So the 
 optimal control problem should have the mixed features of  \emph{ deterministic} and   \emph{  stochastic} controls.
This is indeed the case 
as shown by the main results of this paper, Theorems \ref{T:4.15}, \ref{Necessaryconv} and \ref{T:5.1}.

 \item[3)]
 We establish   a martingale representation theorem for sub-diffusion $B_{L_t}$, which may be of independent interest. It plays a key role in our study of BSDEs  \eqref{y-general} driven by 
sub-diffusions.  

 \item[4)]
When studying the stochastic maximum principle for systems driven by anomalous sub-diffusions, 
we need to consider adjoint  equations \eqref{p} and \eqref{P}, which  are  BSDEs
of the form  \eqref{y-general}  
driven by  the sub-diffusion $B_{L_t}$ with drifts involving both $dt$ and $dL_t$ terms.
We establish  the existence and uniqueness of  
the  solutions for such BSDEs by the martingale representation theorem mentioned in 3) and the contraction mapping theorem. 
 
 \end{enumerate}
  
  \medskip

 \begin{remark} \label{R:1.1} \rm 
 \begin{enumerate}
 \item[(i)] In a recent paper  \cite{NN}, under the assumption that some adjoint BSDEs have solutions, 
  Nane and Ni  studied sufficient stochastic maximum principle
   for SDEs driven by $B_{L_t}$ as well as by a compensated Poisson random measure $\tilde N(dt, dz)$ time-changed by $L_t$, where $L_t$ is the inverse of a purely discontinuous subordinator $S$. 
 Using a Picard's iteration method, they 
  gave an existence and uniqueness result for BSDEs driven $B_{L_t}$ with a ``drift'' or ``generator'' term    $-h(t, L_t, X_{t-}, Z_t) dL_t + \int_{\R\setminus 0} r(t, z) \tilde N (dL_t, dz)$, where $h$ and $r$ are given functions. 
However,  this result is not applicable to the adjoint BSDEs used in their sufficient SMPs \cite[Theorems 3.1 and 4.1]{NN}, where $r(t, z)$ is a part of the solution for the adjoint equations. So the existence of solution to the adjoint BSDEs is a part of the assumptions of the main results in \cite{NN}.

\item[(ii)]  We take this opportunity to point out an issue about the uniqueness of solutions to BSDE 
 in  \cite[Lemma 2.3]{NN}. Since the subordinator $S$ in \cite{NN} does not have a positive drift, $dL_t$ is singular with respect to the Lebesgue measure $dt$ on $[0, \infty)$. Thus from $\int_0^T |u_1(s)-u_2(s)|^2 dL_s =0$, one can not conclude that $u_1(s)=u_2(s)$ for a.e. $s\in (t, T)$ as claimed on \cite[p.200, line 12]{NN}. A similar issue occurs in the definition of Hamiltonian $H$ on \cite[ p.209,  lines 3 and 8]{NN} 
 using  $\frac{dL_t}{dt}$, and, subsequently, in the assumption of $\hat H(x)$ in (4.4) of \cite[Theorem 4.4]{NN} 
 that contains a $\frac{dL_t}{dt}$ term which is 
   required to be   a concave function of $x$ for every $t\in [0, T]$. 
 
 \item[(iii)] While the BSDE considered in \cite[Lemma 2.3]{NN} allows a compensated Poisson Radon measure term but no 
 $dt$ term, ours have a term driven by $dt$.
 In our paper, the existence and uniqueness of the adjoint equations are established, not assumed. We not only give sufficient conditions but also the necessary conditions for a stochastic control to be optimal. 
 \end{enumerate} 
 \end{remark}

  \smallskip
 
In a forthcoming paper \cite{CZ2}, we will study HJB equation for stochastic optimal control for SDEs driven by sub-diffusions. 

\smallskip
 
  The rest of this paper is organized as follows. We formulate the  stochastic control problem in Section \ref{S:2}.  
 In Section \ref{S:3}, overshoot process is introduced  to make the sub-diffusion  Markov.  
  In Section \ref{S:4}, a martingale representation for sub-diffusions is established, which is then used to 
 obtain existence and uniqueness of BSDEs driven by sub-diffusions. This is one of the main results of this paper.
 In Section \ref{S:5},  
  necessary  conditions for optimal control  are
   established 
    using both  
    spiking and convex variational methods.  Section \ref{S:6} is concerned with two sufficient conditions for a control to be optimal, with and without the concavity assumption.
     In Section\ref{S:7}, we apply our main results to  a linear quadratic  system 
      driven by sub-diffusions 
     for which we are able to find the optimal control and 
 its state process explicitly. 
 
     For simplicity, in this paper,  we formulate and state the models and theorems in the setting of one-dimensional
	state space. However, all the results in this paper hold in the high dimensional state spaces by the same argument with some straightforward modifications.  			

\section{The model}\label{S:2} 
 
  Suppose that $B$  is a  standard Brownian motion on $\RR$ starting from $0$,  $S$ is any subordinator that is independent of $B$ with $S_0=0$,
and $L_t:=\inf\{r>0:S_r>t\}$. It will be shown in Theorem \ref{T:3.1} below that
 \begin{equation} \label{e:2.1}
 \wt X_t:=(X_t, R_t):=  \left(x_0+ B_{L_{(t-R_0)^+}},  \, R_0+ S_{L_{(t-R_0)^+}} - t  \right), \quad t\geq 0,
\end{equation}
with $ \wt X_0= (x_0, R_0) \in \RR \times [0, \infty)$   is a time-homegenous Markov process taking values in $\RR\times [0, \infty)$.

 \begin{definition}\label{D:2.1}  \rm
 Let $U$ be a  non-empty Borel subset of $\R$. For each $0\leq s  < T$ and  $a\geq 0$, denote by $\mathcal{U}_{ a} [s, T]$ the set of all 5-tuples  $\left(\Omega,   \mathcal{F} ,    \mathrm P, X, u(\cdot)  \right) $
 satisfying the following conditions:
 \begin{enumerate}
 \item[(i)] $\wt X=\{ \wt X_t; t\in [0, \infty)\}$  is given by  \eqref{e:2.1}  with $\wt X_0=(0, a)$;

  \item[(ii)] $\{u (t, \omega );  t\in [s, T]\}$ defined on $   [s, T] \times \Omega $ is an $\{\sF_{t-s}\}_{t\in [s, T]}$-progressively measurable process 
  taking values in $U$ so that 
   $\E \int_s^T u(t)^2 dt <\infty$, where $\{ \sF_t\}_{t\geq 0} $
  is the minimum augmented filtration generated by $\wt X$.
 \end{enumerate}

Such a   5-tuple   $\left(\Omega,   \mathcal{F} ,    \mathrm P, X, u(\cdot)  \right) $  will be called an admissible control.
We often abbreviate it  as
  $u\in \mathcal{U}_a[s, T]$.     Note that  $\mathcal{U}_a[s, T]$ depends on the open set $U$ but for notational convenience we do
  not include $U$ in its notation. In the following we call $U$ control domain. Observe that $\mathcal{U}_a[s, T]$ is convex if and only if $U$
  is convex. 
  We say $u\in \mathcal{U}_a' [s, T]$ if the filtration $\{\sF_t\}_{t\geq 0}$ in (ii)
  is replaced by the   minimum augmented filtration   $\{\sF'_t\}_{t\geq 0}$  generated by $ X$, the first coordinate process of $\wt X=(X, R)$.
 Clearly, $\sF'_t\subset \sF_t$ for every $t\geq 0$ and so $\mathcal{U}_a' [s, T] \subset \cU_a [s, T]$.
 \end{definition}

\begin{remark}\label{R:2.2}  \rm
\begin{enumerate}
 \item[(i)] For simplicity, we assume the control domain $U$ to be a  non-empty Borel subset of $\R$. It can be replaced by any
separable metric space with the arguments throughout this paper.

\item[(ii)] Note that process $\wt X_t=(X_t, R_t)$ in Definition \ref{D:2.1}
  depends on the initial $a\geq 0 $ of $R_0$, so do the filtrations $\{\sF_t\}_{t\geq 0}$ and $\{\sF'_t\}_{t\geq 0}$. For emphasis, sometimes we denote them
by $\wt X^a_t=(X^a_t, R^a_t)$, $\{\sF^a_t\}_{t\geq 0}$ and $\{ {\sF^{a}_t}'\}_{t\geq 0}$, respectively.  Clearly, $\sF^a_t$ and $ {\sF^{a}_t}' $ are trivial for $t\in [0, a]$.
Hence for each $u\in \cU_a [ s, T]$, $\{u(r); r\in [s, (s+a)\wedge T]\}$ is deterministic.

\item[(iii)]  Note that for $0\leq s<\bar s <T$,  for $u\in \cU_a[s, T]$, its restriction $\{u(r); r\in [\bar s, T]\}$ on the time interval $[\bar s, T]$
is in general not a member in $\cU_a [\bar s, T]$.   \qed
 \end{enumerate}
\end{remark}
\medskip

Fix $a\geq 0$. Given $u\in \mathcal{U}_a [s, T]$  and $x_0\in \R$, consider the following SDE driven by the anomalous sub-diffusion $X_t$ for
$x^u=x^{u, s, x_0, a}$:
 \begin{equation}\label{e:state}
\left\{\begin{aligned}
dx^u (t)=&\ b(t, x^u (t), u(t))dt  +\sigma(t, x^u (t), u(t) )dB_{L_{(t-s-a)^+}}   \quad \hbox{for } t\in[s,T],    \\
 x^u (s)=&\ x_0.
\end{aligned}
\right.
\end{equation}

\medskip

When $b$ and $\sigma$ are Lipschitz in $x$ uniform in $(t, u)$, it is well known that the SDE \eqref{e:state} has a unique strong as well as weak solution for every
$u\in \mathcal{U}_a [s, T]$  and $x_0\in \R$. The strong solution $\{ x^u (t); t\in [s, T]\}$ is a continuous process that is progressively
measurable with respective to the filtration  $\{\sF_{t-s}; t\in [s, T]\}$. When $u\in \cU'_a [s, T]$, then  the strong solution $\{x^u (t); t\in [s, T] \}$ is a continuous process that is progressively
measurable with respective to the smaller filtration  $\{\sF'_{t-s}; t\in [s, T]\}$.
  Many times, it is  natural to consider controls from  $\mathcal{U}_a' [s, T] $ as the driving force for \eqref{e:state} is the sub-diffusion $X_t=X_0+B_{L_{(t-a)^+}}$ and one does not have  the complete information of $\wt X=(X, R)$. In the remaining part of this paper, we mainly consider 
optimal controls over $u$ in $\mathcal{U}_a' [s, T] $.

 \begin{remark}\rm
 \begin{enumerate}
 \item[(i)]  If $a\geq T-s$, then $(t-s-a)^+=0$ for every $t\in [s, T]$. In this case, the SDE \eqref{e:state} degenerates into a deterministic ODE.
  For stochastic optimal controls considered in this paper, the non-trivial case is when $a\in [0, T-s)$ though we impose this restriction  in the paper.

 \item[(ii)]
Note that    $ L_{(t-a)^+} = \<X  \>_t$   is the quadratic variational process of the continuous
 local martingale $X$, $L_{(t-a)^+}$ is $\sF'_t$-measurable.  By \cite[Lemma 0.4.8]{RY}
 $$
S_t=\inf\{r>0: L_r>t\} = \inf\{r>0: L_{(r -a)^+} >t\}-a  \quad \hbox{for every } t\geq 0.
 $$
 Hence each $S_t+a$ is an $\{\sF'_s\}_{s\geq 0}$-stopping time and
 $ R_t:=   S_{L_{(t-a)^+}}+a - t$
 is  not $\sF'_t$-measurable as its valued depends on the  future value of  $L_s$  beyond $L_{(t-a)^+}$.  This shows that the filtration $\{\sF'_t\}_{t\geq 0}$ is a proper sub-filtration of $\{\sF_t\}_{t\geq 0}$.
 For some part of the theory developed in this paper, it works for general admissible control $u\in \cU_a [s, T]$. However our martingale representation theorem
 for  subdiffusion $X$ requires the terminal   random variables  be in $\sF'_T$; see Theorem \ref{T:3.2}. So any results that use Theorem \ref{T:3.2} will require the admissible control
 $u$ from the  smaller class $\cU_a' [s, T]$.      
  \item[(iii)]  Our formulation of the state process $x^u(t)$ of form \eqref{e:state} is motivated by the stock price example in bear market mentioned in the Introduction. 
 We could allow an addition $d L_{(t-s-a)^+}$ term for the state process $x^u(t)$ in \eqref{e:state}, that is, $x^u (t)$ is governed by
 \begin{eqnarray}\label{e:state2}
 dx^u (t) &=& b_1(t, x^u (t), u(t))dt +b_2(t, x^u (t), u(t))dL_{(t-s-a)^+}   \\
&&   +\sigma(t, x^u (t), u(t) )dB_{L_{(t-s-a)^+}}   \quad \hbox{for } t\in[s,T]   \nonumber
  \end{eqnarray}
  with $x^u (s)=  x_0$. All the arguments and the results in this paper carry over with appropriate modifications.
  Indeed, in    
 \cite{ZhangChen-FB}, we have considered the case where the state process is in such a generality. 
 For simplicity, we choose not to include such an extension having an extra term $dL_{(t-s-a)^+}$ for the state process in this paper. 
 However, for the adjoint backward equations \eqref{p} and \eqref{P}, it is important to have the $dL_{(t-s-a)^+}$ term. 
 
 \item[(iv)] In this paper, the subordinator $S$ is assumed to have a positive drift $\kappa$ for two reasons. 
 First, when the L\'evy measure $\nu$ of $S$ is null, the sub-diffusion $B_{L_t}$ reduces to Brownian motion
 and hence our esults in particular cover the classical models driven by Brownian motion as treated in literature, see, \cite{YZ}, for instance. 
 The main results of this paper can also be viewed as a stability result in the sense that if the L\'evy measure $\nu$ tends to zero, the corresponding  SMP converges to that for the Brownian model.
 Secondly, when $\kappa >0$, $dL_t$ is absolutely continuous with respect to the Lebesgue measure $dt$ as mentioned in \eqref{e:1.1}.
 This simplifies the formulation of the SMPs in Theorems \ref{T:4.15} and \ref{Necessaryconv}.
When $\kappa =0$, it might still be possible to obtain suitable SMP for sub-diffusions, but the formulation will be more involved
which needs  to deal with the  random measure $dL_t$ on $[0, T]$ 
that  is singular with respect to the Lebesgue measure $dt$ on $[0, T]$;
see the proof of Theorem \ref{T:4.15} for an example.  For this reason and for the sake of presenting our appraoch for stochastic optimal control for systems driven by anomalous sub-diffusions as transparent as possible,
we choose to assume $S$ having a positive drift. 
 \end{enumerate}
 \end{remark}

   Consider the following cost functional  for control $u\in \cU_a [s, T]$:
   \begin{eqnarray}\label{e:2.3}
J(s, x_0, u, a  )&=& \EE \bigg[\int_{s}^{T}  f\left  (t, x^{u, s, x_0,  a}  (t),
u  (t) \right)   dt+h(x^u (T)) \bigg].
\end{eqnarray}
 The optimal control problem is to find the control $u^\ast \in \cU'_a [s, T]$  (respectively,  $u^\ast \in \cU_a [s, T]$) to minimize the above cost functional
  \begin{equation}\label{JSSC-Eq5}
 J(s, x_0, u^\ast,  a  )  =\inf_{u \in \mathcal{U}'_a[s, T]} J(s, x_0, u, a  )
\end{equation}
(respectively,   $J(s, x_0, u^\ast,  a  )  =\inf_{u \in \mathcal{U}_a[s, T]} J(s, x_0, u, a  )$).
  Unless otherwise specified, in the remaining part of this paper, we take the initiala time $s$ to be zero. 
  
 \section{Overshoot process}\label{S:3}

  Note that for each fixed $t>0$,  $S_{L_t}>t$ happens with positive probability.
 On $\{S_{L_t}>t\}$, the inverse local time $L_s$ and, consequently, 
 the sub-diffusion $B_{L_s}$ remain flat during the time interval $[t, S_{L_t} ]$. 
 We call $S_{L_t}-t$ an overshoot process with initial value 0. It measures how much time it would take for the anomalous sub-diffusion $B_{L_s}$ to wake up   from time $t$. 
The anomalous sub-diffusion 
$B_{L_{(t-a)^+}}$ itself  is not Markov.  As the next theorem shows,   
 we can add the overshoot process  to make it Markov.   
This Markov property will be used in the proof of Proposition \ref{P:3.2} on a property of inverse local time $L_t$.
It also plays an important role in our study \cite{CZ2} of the 
dynamic programming principle and the Hamilton–Jacobi–Bellman equations for systems driven by anomalous sub-diffusions.

 \begin{theorem}\label{T:3.1} 
  Suppose that $B$  is a  standard Brownian motion on $\RR$ starting from $0$,  $S$ is any subordinator that is independent of $B$ with $S_0=0$,
and $L_t:=\inf\{r>0:S_r>t\}$.
Then
\begin{equation} \label{e:5.1}
\wt X_t:=(X_t, \, R_t):= \left(x_0+ B_{L_{(t-R_0)^+}},  \, R_0+ S_{L_{(t-R_0)^+}} - t  \right), \quad t\geq 0,
\end{equation}
with $\wt X_0= (x_0, R_0) \in \RR \times [0, \infty)$   is a time-homegenous Markov process taking values in $\RR\times [0, \infty)$.
\end{theorem}

\noindent{\bf Proof.} Note that $(B, S)$ is a L\'evy process. Thus for any $t\geq 0$,
 \begin{equation}\label{e:2.2} 
\{\wt X_{t+s}- \wt X_t; s\geq 0\}= \left\{ \left(B_{L_{(t+s-R_0)^+} - L_{(t-R_0)^+}}, S_{L_{(t+s-R_0)^+}- L_{(t-R_0)^+} }-s \right) \circ \theta_{L_{(t-R_0)^+} } ; s\geq 0 \right\},
\end{equation} 
where  $\theta_r$ is  the time shift operator for the L\'evy process  $(B, S)$.

Suppose $\wt X_t=  (X_t, R_t)= \left(x_0 + B_{L_{(t-R_0)^+}},  S_{L_{(t-R_0)^+}}+R_0 -t \right)=(x, a)$ and $s\geq 0$.
If $t+s \leq R_0$, then $a=R_0-t$ and
$$ L_{(t+s-R_0)^+} - L_{(t-R_0)^+} = 0  =L_{(s-a)^+} \circ \theta_{L_{(t-R_0)^+} }  ; $$
while if $ {t+s} > R_0 $,
\begin{eqnarray*}
L_{(t+s-R_0)^+} - L_{(t-R_0)^+} &=& \inf\{r>0: S_r>t+s -R_0 \}  - L_{(t-R_0)^+}   \\
&=&  \inf \left\{r>0: S_{r+   L_{(t-R_0)^+}  } -S_{L_{(t-R_0)^+} } >s +t -R_0-S_{L_{(t-R_0)^+} } \right\} \\
&=&  \inf \left\{r>0: S_r\circ \theta_{L_{(t-R_0)^+} } >s -a  \right\} =L_{(s-a)^+} \circ \theta_{L_{(t-R_0)^+} }.
\end{eqnarray*}
 Thus we have for any $t, s>0$, 
 \begin{equation}
L_{(t+s-R_0)^+} - L_{(t-R_0)^+}  =L_{(s-R_t)^+} \circ \theta_{L_{(t-R_0)^+} }.
\end{equation}
  This together with \eqref{e:2.2} implies that   for each fixed $t\geq 0$,
 \begin{equation}\label{e:5.2a}
\{\wt X_{t+s}- \wt X_t; s\geq 0\}= \left\{ \left(B_{L_{( s-R_0)^+}} , S_{ L_{(s-R_0)^+}} -s  \right) \circ \theta_{L_{(t-R_0)^+} } ; s\geq 0 \right\},
\end{equation}
This shows that   the conditional distribution of $\{\wt X_{t+s}; s\geq 0\}$ given $\sF_t$ has the same distribution
as $\{\wt X_s; s\geq 0\}$ starting from the random position $X_t=(x, a)$ at time $s=0$.  \qed

\medskip

Now suppose $S=\{S_t; t\geq 0\}$ is a subordinator with   drift $\kappa >0$ and L\'evy measure $\nu$ starting from 0.  
  Denote by $U(dx)$ the potential measure of the subordinator; that is, for any $f\geq 0$ on $[0, \infty)$,
  $$
  \int_{[0, \infty)}  f(x) U(dx) = \E \int_0^\infty f(S_t) dt.
  $$
  Clearly, for each $x>0$, $\E \left[  L_x \right]  = U([0, x]) =:U(x)$. In the literature, $U(dx)$ and $U(x)$ are also called   the
  renewal measure and  the renewal function, respectively.
  Since $S$ has positive drift $\kappa >0$, by a result due to Reveu (see \cite[Proposition 1.7]{Ber2}),
  $U(dx)$ is absolutely continuous with respect to the Lebesgue measure on $[0, \infty)$,
  has a strictly positive continuous density function ${\vartheta} (x)$ on $[0, \infty)$ with ${\vartheta} (0)=1/\kappa$,
  and
\begin{equation}\label{e:3.1a}
  {\vartheta}(x) =  \kappa^{-1}  \, \PP (\hbox{there is some $t\geq 0$ so that } S_t=x)  \quad \hbox{for }  x\geq 0.
 \end{equation}
 Furthermore,  one has  (see, e.g., \cite[Theorem 5]{Ber1}),
\begin{equation}\label{e:3.1}
  \PP (S_{L_x} =x) = \kappa \,  {\vartheta}(x)  \quad \hbox{for  every }  x > 0.
  \end{equation}
  Consequently,   we have by the bounded convergence theorem that for any $t\geq 0$,
   \begin{equation} \label{e:3.2a}
  \lim_{s\to 0+} \frac{\E \left[ L_{t+s}-L_t\right]}{s}
  =  \lim_{s\to 0+}  \frac{U(t+s)-U(t)}{s}=  \lim_{s\to 0}  \frac{\int_t^{t+s} {\vartheta}(x) dx} {s}={\vartheta}(t)  .
    \end{equation}
  In particular,
  \begin{equation} \label{e:3.2}
  \lim_{t\to 0} \frac{\E \left[ L_t\right]}{t}
   ={\vartheta} (0) = 1/\kappa.
    \end{equation}
 
\medskip

\begin{prop}\label{P:3.2} 
Suppose $S=\{S_t; t\geq 0\}$ is a subordinator with  drift $\kappa >0$ starting from $0$.  
Then for any $R_0\geq 0$ and $t\geq 0$, 
\begin{equation} \label{e:2.9} 
 \lim_{s\to 0+} \frac 1s  \E \left[ L_{(t+s-R_0)^+} - L_{(t-R_0)^+}   \big| \sF_t \right] 
 = \kappa^{-1} \1_{\{R_t=0\}}   = \frac{dL_{(t-R_0)^+}}{dt}, 
\end{equation} 
where  $\{ \sF_t\}_{t\geq 0} $
  is the minimum augmented filtration generated by $\wt X_t= (X_t, R_t)$.
 \end{prop}

\pf  Note that \eqref{e:5.2a} in particular implies that for any $t, s>0$, 
$$
\E \left[ L_{(t+s-R_0)^+} - L_{(t-R_0)^+}   | \sF_t \right] = \left( \E {L_{(s-a)^+}} \right)  \big|_{a=R_t}.
$$
 Thus by \eqref{e:3.2} we have for every $t\geq 0$, 
$$  \lim_{s\to 0+} \frac 1s  \E \left[ L_{(t+s-R_0)^+} - L_{(t-R_0)^+}   \big| \sF_t \right] 
= \lim_{s\to 0+} \frac 1s  \left( \E {L_{(s-a)^+}} \right)  \big|_{a=R_t} = \kappa^{-1} \1_{\{R_t=0\}}.
$$
By \eqref{e:1.1}, $\frac{dL_{(t-a)^+}}{dt} = K_{(t-a)^+}$ exists a.s. with $0\leq K_{(t-a)^+}\leq 1/\kappa$. 
It follows from the dominated convergence theorem that 
$$
K_{(t-R_0)^+}= \lim_{s\to 0+} \frac 1s  \E \left[ L_{(t+s-R_0)^+} - L_{(t-R_0)^+}   \big| \sF_t \right] 
= \kappa^{-1} \1_{\{R_t=0\}}.
$$
This establishes the proposition.

\qed 
 
\section{Martingale representation theorem and BSDEs} \label{S:4}

Let $a\geq 0$.   To prepare for the SMP,  we need to establish the existence and uniqueness for the following  type
backward stochastic differential equation driven by
$B_{L_t}$ on $[0, T]$ for any $T>0$:
\begin{equation}\label{y-general}
 d Y_t   =h_1 (t,Y_t )dt  + h_2 (t,Y_t,Z_t)dL_{(t-a)^+} +Z_t d B_{L_{(t-a)^+}}  \quad \hbox{with }
Y(T)  = \xi,
 \end{equation}
 where $\xi \in L^2(\sF'_T) $, the space of square integrable  $\sF'_T$-measurable random variables.
 (Recall that $\sF'$  is the augmented filtration generated generated by $ X$.)

First, we need the following  integral representation of square integrable random variables with respect to  the subdiffusion  $B_{L_{(t-a)^+}}$.
The following result holds for any subordinator $S$; that is, we do not need to assume that $S$ has a positive drift $\kappa>0$.

\begin{theorem}\label{T:3.2}
 For each  $a\geq 0$, $T\in (0, \infty]$ and $\xi \in L^2(\sF'_T)$, there exists an  
   $\{\sF'_t\}_{t\in [0, T]}$-progressively measurable 
  process $H_s$   with   $\E \int_0^T H_s^2 dL_{(s-a)^+} <\infty$ so that
\begin{equation}\label{e:3.12}
\xi = \E [\xi] + \int_0^T H_s d B_{   L_{(s-a)^+ }   }.
\end{equation}
Such $H$ is unique in the sense that if $H'$ is another 
 $\{\sF'_t\}_{t\in [0, T]}$-progressively measurable   process,
 then $\E \int_0^T |H_s-H'_s|^2 dL_s =0$.
\end{theorem}

\pf  
By \cite{CZ},  for each $t>  0$ and $\lambda >0$,
$ \E \left[  e^{ \lambda L_t } \right] <\infty$. For any $n\geq 1$, $0\leq t_0< t_1 <t_2 <\cdots <  t_n \leq T$ and $\lambda_i \in \R$ for $1\leq i \leq n$,
\begin{equation} \label{e:3.13}
Z_t := \exp \left( \int_0^t f(s) d B_{L_{(s-a)^+}}  -\frac12 \int_0^t f^2  (s) d L_  {(s-a)^+   }   \right)
\end{equation}
is a square integrable $\{\sF'_t\}_{t\in [0, T]}$-martingale, where
  $f(t) := \sum_{i=1}^n  \lambda_i \1_{(t_{i-1} , t_i]} (t)$.
  Note that $Z_t$ satisfies
  \begin{equation} \label{e:3.14}
  Z_t =1 + \int_0^t Z_s f(s) dB_ {   L_{(s-a)^+ }    }    .
 \end{equation}
     By the same proof as that for \cite[Lemma V.3.1]{RY},  random variables $Z_T$ of this kind is dense in $L^2(\sF'_T)$.

Let $\sH_T$ be the subspace of elements $\xi$ in $L^2(\sF'_T)$ that admits representation \eqref{e:3.12}.
For $\xi \in \sH_T$, we have
$$
\E [\xi^2 ] = (\E [ \xi ] )^2 + \E \left[ \int_0^T H_s^2 dL_ {(s-a)^+ }    \right] .
$$
 If $\{\xi_n; n\geq 1\}$ is a Cauchy sequence in $\sH_T$, its   corresponding 
   $\{\sF'_t\}_{t\in [0, T]}$-progressively measurable 
 sequence $\{H^{(n)}_s; n\geq 1\}$
satisfies
$$
 \lim_{n, k\to \infty} \E \left[ \int_0^T (H^{(n)}_s-H^{(k)}_s)^2 dL_ {(s-a)^+ }   \right] =0 .
 $$
 Then there is 
   an $\{\sF'_t\}_{t\in [0, T]}$-progressively measurable 
 process  $H_s$ with   $\E \int_0^T H_s^2 dL_s <\infty$
 so that
 $$
 \lim_{k\to \infty} \E \left[ \int_0^T (H^{(k)}_s-H_s)^2 dL_ {(s-a)^+ }    \right] =0 ;
 $$
cf. \cite[p.130]{RY}.
 Consequently,
 $$
 \lim_{k\to \infty} \E \left[   \left( \int_0^T H^{(k)}_s d B_{L_ {(s-a)^+ }    }  -   \int_0^T H_s d B_{L_ {(s-a)^+ }   } \right)^2 \right]
 = \lim_{k\to \infty} \E \left[ \int_0^T (H^{(k)}_s-H_s)^2 dL_ {(s-a)^+ }    \right] =0 .
 $$
 It follows that the $L^2$-limit $\xi$ of $\{\xi_n; n\geq 1\}$ admits the representation
 $  \xi = \E [\xi] + \int_0^T H_s d B_{L_{(s-a)^+ }    }$, which shows that $\xi \in \sH_T$;
 that is, $\sH_T$ is closed with respect to the $L^2$-norm.
 Since the  random variables $Z_T$ of the form \eqref{e:3.13} are dense in $L^2(\sF'_T)$ and, in view of \eqref{e:3.14},
 are members of  $\sH_T$, we conclude that $L^2(\sF'_T)=\sH_T$. This proves the first part of the theorem.

 Suppose that $H'$ is another   
   $\{\sF'_t\}_{t\in [0, T]}$-progressively measurable process
   with   $\E \int_0^T (H'_s)^2 dL_ {(s-a)^+ }     <\infty$ so that
 $ \xi = \E [\xi] + \int_0^T H_s d B_{L_ {(s-a)^+ }    } $.
 Then $\int_0^T (H_s -H'_s) d B_{L_  {(s-a)^+ }    } =0$. Consequently,
 $$
 \E \int_0^T |H_s-H'_s|^2 dL_ {(s-a)^+ }     = \E \left[ \left( \int_0^T (H_s -H'_s) d B_{L_ {(s-a)^+ }   } \right)^2 \right] =0.
 $$
 This completes the proof of the theorem.
 \qed

\medskip

For  $a\geq 0$ and  $\beta >0 $,  we can define a Banach norm 
$\| \cdot  \|_{ \mathcal M_{a, \beta} [0,T]}$  on the space  
\begin{eqnarray*}
\mathcal M [0,T] 
 &=& \Big\{  (Y, Z): \hbox{ $Y$ and $Z$ are $\{\sF'_t\}_{t\in [0, T]}$-progressively measurable processes on $[0, T]$ with }  \\
&& \hskip 0.8truein 
 \EE  \Big[    \int_0^T  | Y(t)|^2dt  +  \int_0^T   | Z(t)|^2  dL_ {(t-a)^+ }  \Big] <\infty \Big\} 
\end{eqnarray*}
by 
\begin{eqnarray}
\| (  Y(\cdot), Z(\cdot)    )  \|_{ \mathcal M_{a, \beta} [0,T]} \triangleq
\left(    \EE  \left[    \int_0^T e^{2\beta t}   | Y(t)|^2dt  +  \int_0^T   | Z(t)|^2 e^{ 2\beta t}  dL_ {(t-a)^+ }  \right]  \right)^{1/2}.
\end{eqnarray}
Clearly for finite $T>0$,   all these norms are equivalent on $\mathcal M [0,T] $. However, later we like to choose $\beta$ to be sufficiently large so that 
certain map becomes a contraction map under the norm $\| \cdot  \|_{ \mathcal M_\beta [0,T]}$;
see the proof of Theorem \ref{T:BSDE}. When   $\mathcal M [0,T]$ is equipped with the norm $\| \cdot      \|_{ \mathcal M_\beta [0,T]}$,
we may denote $\mathcal M [0,T]$ by $\mathcal M_{a, \beta} [0,T]$ for emphasis. For notational simplicity, we will denote $\| \cdot      \|_{ \mathcal M_\beta [0,T]}$
by $\| \cdot      \|_{ \mathcal M_\beta}$ when there is no danger of ambiguity.

 \medskip

  \begin{hypothesis}\label{HPeuy}
For any $(y,z)\in\RR\times\RR $, $h_1(t,y ) $  and  $ h_2(t,y,z ) $  are 
 $\{ \sF'_t \}_{t\geq 0}$-progressively measurable random processes 
with
 $\E \int_0^T \left( |h_1 (s, 0 )|^2 + |h_2 (s, 0,0)|^2 \right) ds <\infty$.   Moreover, there exists a constant $C > 0$
so that $\PP$-a.s., 
$$
|h_1(t,x_1  )-h_1(t,x_2 )|\leq C |x_1 -x_2| \ \hbox{ and } \
|h_2(t,x_1 ,z_1 )-h_2(t,x_2,  z_2 )|\leq C\( |x_1 -x_2|+|  z_1 - z_2|\)
$$
  for  any $ t \in [0,T]$ and $x_i, z_i \in \RR$ with $i=1, 2$. 
\end{hypothesis}

\medskip

\begin{definition}\label{YZsolution}
A pair of process $ (Y(\cdot), Z(\cdot)) \in \mathcal M[0,T]$ is called an adapted solution of (\ref{y-general}) if  for any $t\in [0, T]$,
\begin{equation} \label{e:4.6}
Y_t   =\xi-\int_t^T   h_1 (s,Y_s)ds-\int_t^T   h_2(s,Y_s, Z_s)dL_ {(s-a)^+ }  -\int_t^T  Z_s d B_{L_ {(s-a)^+ }   }   \quad  \PP  \hbox{-a.e.} 
\end{equation}
 \end{definition}

\begin{lemma}\label{L:4.4}
Suppose that $ (Y(\cdot), Z(\cdot)) \in \mathcal M[0,T]$ is an adapted solution of (\ref{y-general}). Then
there is a positive constant $c$ that depends on $\kappa$, $T$ and the Lipschitz constant $C$ in Hypothesis \ref{HPeuy}
so that  
\begin{eqnarray}\label{e:4.7}
\EE \left[ \sup_{t\in [0, T]} |Y_t|^2 \right] 
&\leq&   c  \E \left[ |\xi|^2 + \int_0^T \left( |h_1(s, 0)|^2  + |h_2(s, 0, 0)|^2 \right)d s \right] \nonumber \\
&& + c\, \E \int_0^T |Y_s|^2 ds +   c\, \E \int_0^T |Z_s|^2 dL_{(s-a)^+}.
\end{eqnarray}
\end{lemma}

\pf It follows from \eqref{e:4.6}, Hypothesis \ref{HPeuy},  Doob's maximal inequality and \eqref{e:1.1}, that
\begin{eqnarray*}
&& \EE \left[ \sup_{t\in [0, T]} |Y_t|^2 \right]  \\
&\leq &   \E \left[\left( |\xi| + \int_0^T |h_1(s, Y_s)| ds + \int_0^T |h_2(s, Y_s, Z_s)| dL_{(s-a)^+}
+ 2\sup_{t\in [0, T]} \Big|\int_0^t Z_s dB_{L_{(s-a)^+}} \Big| \right)^2 \right] \\
&\leq & 4 \E \left[ |\xi|^2 \right] + 4 \E \left[ T \int_0^T |h_1(s, Y_s)|^2 ds \right]
+ 4 \E \left[L_{(T-a)^+}  \int_0^T |h_2(s, Y_s, Z_s)|^2 dL_{(s-a)^+} \right] \\
&& +64 \E \left[ \int_0^T |Z_s|^2 dL_{(s-a)^+} \right] \\ 
&\leq & 4 \E \left[ |\xi|^2 \right] + 8TC^2 \E \left[ \int_0^T \left( |h_1(s, 0)|^2 + |Y_s|^2 \right) ds \right] \\
&& + 12 \kappa^{-1}T C^2 \E \left[ \int_0^T \left( |h_2(s, 0, 0)|^2 + |Y_s|^2 +|Z_s|^2\right)  dL_{(s-a)^+} \right]  
 +64 \E \left[ \int_0^T |Z_s|^2 dL_{(s-a)^+} \right] \\
&\leq & c \E \left[ |\xi|^2 + \int_0^T \left( |h_1(s, 0)|^2  + |h_2(s, 0, 0)|^2 \right)d s \right] 
  +c\E \int_0^T |Y_s|^2 ds + c \E \int_0^T |Z_s|^2 dL_{(s-a)^+}. 
\end{eqnarray*}
\qed

\begin{theorem}\label{T:BSDE}
Suppose Hypothesis \ref{HPeuy} holds. Then for any given $\xi \in  L^2 (\sF'_T)$,
  the BSDE  \eqref{y-general}  admits an adapted solution $(Y(\cdot), Z(\cdot)) \in \mathcal M[0,T]$.
  The solution is unique in the sense that if $(\wt Y(\cdot), \wt Z(\cdot)) \in \mathcal M[0,T]$ is another solution of \eqref{y-general},
  then $\wt Y_t=Y_t$ for all $t\in [0, T]$ with probability one and $\E \int_0^T |Z_s -\wt Z_s|^2 dL_ {(s-a)^+ }    =0$.
\end{theorem}

\pf   
  Let  $\xi \in L^2(\sF'_T)$.
   Given $(y_t, z_t)\in  \mathcal M[0,T]$,
 consider the following BSDE:
\begin{eqnarray}\label{e:3.19}
 d Y_t   =&h_1(t, y_t )dt + h_2(t, y_t, z_t)d L_ {(t-a)^+ }   +Z_t d B_{L_ {(t-a)^+ }  }   \quad \hbox{with }  \ 
Y_T  =   \xi   .
\end{eqnarray}
 Define 
 $$
 \eta = \xi- \int_0^T h_1(t, y_t )dt -\int_0^T  h_2(t, y_t, z_t)dL_ {(t-a)^+ } .
 $$
Note that under Hypothesis  \ref{HPeuy} and by \eqref{e:1.1}, 
\begin{eqnarray*}
\E [\eta^2]
 &\leq& 3\E [ \xi^2] + 3T \E \int_0^T |h_1(t, y_t)|^2 dt + 3\E \left[\left( \int_0^T |h_2(t, y_t, z_t) | dL_{(t-a)^+} \right)^2 \right] \\
&\leq & 3\E [ \xi^2] + 3T  C^2 \E \int_0^T (|h_1(t, 0)|+|y_t|)^2 dt \\
&& + 3C^2 \E \left[\left( \int_0^T ( |h_2(t, 0, 0)| + |y_t| +|z_t|) dL_{(t-a)^+} \right)^2 \right]  \\
&\leq &  3\E [ \xi^2] + cT\E \int_0^T |h_1 (t, 0)|^2 dt +  cT \kappa^{-2}\E \int_0^T |h_2 (t, 0, 0)|^2 dt  \\
&& + cT (1+T\kappa^{-1}) \E \int_0^T |y_t|^2 dt +  c   T  \kappa^{-1} \E \left[ \int_0^T |z_t|^2 dL_{(t-a)^+}  \right] 
< \infty.
\end{eqnarray*}
 
By Theorem \ref{T:3.2} there exists  
an $\{\sF'_t\}_{t\in [0, T]}$-progressively measurable 
process $Z$
 with   $\E \int_0^T  Z_s^2 d L_{(s-a)^+ }  <\infty$ so that
 $\eta = \E [\eta] + \int_0^T Z_s d B_{L_{(s-a)^+ }    }$.

Define 
$$ Y_t= \E \eta  + \int_0^t h_1(s, y_s)ds + \int_0^t  h_2(s, y_s, z_s)dL_ {(s-a)^+ }  +   \int_0^t  Z_s dB_{L_ {(s-a)^+ }   }.
$$
Then for $t\in [0, T]$, 
\begin{eqnarray*}
 Y_t &=& \eta -  \int_0^T Z_s d B_{L_{(s-a)^+ }    }  + 
 \int_0^t h_1(s, y_s )ds + \int_0^t  h_2(s, y_s, z_s)dL_ {(s-a)^+ }  +   \int_0^t  Z_s dB_{L_ {(s-a)^+ }   } \\
 &=&  \xi  -  \int_t^T h_1(s, y_s )ds - \int_t^T  h_2(s, y_s, z_s)dL_ {(s-a)^+ } -  \int_t^T  Z_s dB_{L_ {(s-a)^+ }   } .
 \end{eqnarray*}
 Thus  $(Y_t, Z_t) \in  \mathcal M[0,T]$
solves BSDE  \eqref{e:3.19}. 
Suppose that $(Y'_t, Z'_t) \in  \mathcal M[0,T]$ is another solution of BSDE  \eqref{e:3.19}.
Then $d(Y_t-Y'_t)= (Z_t-Z'_t) dB_{L_ {(t-a)^+ }   }$ with $Y_T-Y'_T=0$. 
  It follows from  Theorem \ref{T:3.2} that $\E \int_0^T(Z_s-Z_s')^2 d L_ {(s-a)^+ } =0$ and 
  consequently $Y_t=Y_t'$ for all $t\in [0, T]$. 
  This shows that BSDE  \eqref{e:3.19}  has a unique solution in $\mathcal M[0,T]$.

 The above defines a map    $\Phi:  \mathcal M [0,T]\rightarrow  \mathcal M  [0,T] $ by sending $(y,z)$ to the solution $ (Y,Z)$
 of \eqref{e:3.19}.   We next show that the map $\Phi$ is contractive with respect to the Banach norm $\| \cdot \|_{ \mathcal  M_{a, \beta} [0,T]}$
 for sufficiently large $\beta >0$.  
 For $(y,z)$ and $(\wt y, \wt z)\in  \mathcal M_\be [0,T]$, let $(Y, Z)= \Phi (y, z)$ and $(\wt Y, \wt Z) = \Phi ( \wt y, \wt z)$.
 For notational simplification, let
 $$
 \wh Y_t:= Y_t - \wt Y_t, \quad \wh Z_t:= Z_t - \wt Z_t, \quad   \wh y_t:= y_t - \wt y_t, \quad \wh z_t:=z_t - \wt z_t,
 $$
 and
 $$
 \wh h_1 (t):=  h_1(t,   y_t)- h_1(t, \wt y_t )  ,  \qquad  \wh h_2 (t):=  h_2(t,   y_t,   z_t)- h_2(t, \wt y_t, \wt z_t)  .
 $$
   By \eqref{e:3.19},
 \begin{eqnarray*}
 \wh Y_t &=&    \int_t^T (   h_1(s, \wt y_s ) - h_1(s,   y_s ) ) ds + \int_t^T (   h_2(s, \wt y_s, \wt z_s) - h_2(s,   y_s,   z_s) ) dL_ {(s-a)^+ }    + \int_t^T  (\wt Z_s -Z_s) dB_{L_ {(s-a)^+ }  } \\
 &=& - \int_t^T \wh h_1 (s) ds - \int_t^T \wh h_2 (s) dL_{(s-a)^+ }    - \int_t^T \wh Z_s dB_{L_ {(s-a)^+ }    } .
\end{eqnarray*}
 Let $\beta>0$, whose value will be taken to be sufficiently large later. 
  Since $\wh Y_T =0$,  applying Ito's formula to $|\wh Y_s  |^2 e^{2 \beta s}$ and evaluating  at $s=t$ and $s=T$ yields
  \begin{eqnarray*}
 && |\wh Y_t|^2 e^{2\beta t} + \int_t^T |\wh Z_s|^2 e^{2\beta s} dL_ {(s-a)^+ }    \\
 &=& -  2 \int_t^T \left(  \beta |\wh Y_s|^2 +\<\wh Y_s, \wh h_1 (s)\> \right)e^{2\beta s} ds -  2 \int_t^T    \<\wh Y_s, \wh h_2 (s)\>   e^{2\beta s} dL_  {(s-a)^+ }
 - 2 \int_t^T \<\wh Y_s, \wh Z_s\> e^{2\beta s} dB_{L_ {(s-a)^+ }    } .
   \end{eqnarray*}
    Denote the continuous local martingale 
$t\mapsto \int_0^t  \<\wh Y_s, \wh Z_s\> e^{2\beta s} dB_{L_ {(s-a)^+ } }$ by $M_t$, and 
its  quadratic variational process by  $\< M\>$.
From \eqref{e:4.7},   
  \begin{eqnarray*}
\E \left[ \< M\>_T^{1/2} \right]
&=&  \E  \left[ \left( \int_0^T  \left| \<\wh Y_s, \wh Z_s\> e^{2\beta s} \right|^2 d  L_ {(s-a)^+  }\right)^{1/2} \right] \\
   &\leq & e^{2\beta T} \E \left[ \left( \sup_{t\in [0, T] }|\wh Y_t| \right) 
	\left(\int_0^T   | \wh Z_s|^2  d  L_ {(s-a)^+}\right)^{1/2} \right] \\
	&\leq & \frac12 e^{2\beta T} \E \left[   \sup_{t\in [0, T] }| \wh Y_t|^2 +
	 \int_0^T   |\wh Z_s|^2  d  L_ {(s-a)^+} \right] 
	< \infty.
  \end{eqnarray*} 
	Thus by the Burkholder-Davis-Gundy inequality, 
	\begin{equation}\label{e:4.9}
	\E \left[ \sup_{t\in [0, T]}|M_t|\right]\leq C \E \left[ \<M\>_T^{1/2}\right] <\infty . 
	\end{equation}
	Hence $\{M_t; t\in [0, M]\}$ is a uniformly integrable martingale.
 
Under Hypothesis \ref{HPeuy},
 $|\wh h_1(s)| \leq C |\wh y_s|$ and $|\wh h_2 (s)| \leq C (  |\wh y_s| +|\wh z_s|)$.
 Noting also that $dL_s \leq \kappa^{-1} ds$,
   we can follow the proof for Theorem 3.2 on p.356-358 in \cite{YZ} to establish that $\Phi$ is a contraction map under $\|  \cdot  \|_{ \mathcal M_{a, \beta} [0,T]}$
     by choosing $\beta >0$ large enough. From the last display, we have 
\begin{eqnarray}\label{esti-norm1}
&& |\wh Y_t|^2 e^{2\beta t} + \int_t^T |\wh Z_s|^2 e^{2\beta s} dL_ {(s-a)^+ }   \nonumber  \\
  &\le&\int^T_t\left(-2\be\widehat{Y}^2_s+2|\wh{Y}_s|C\left(|\wh{y}_s|\right)\right)e^{2\be s}ds
+\int^T_t\left( 2|\wh{Y}_s|C\left(|\wh{y}_s|+|\wh{z}_s|\right)\right)e^{2\be s} dL_  {(s-a)^+ }
\nonumber\\
&& -\int^T_t2e^{2\be s}\wh{Y}_s\wh{Z}_sdB_{L_ {(s-a)^+ } }\nonumber\\
&\le&\int^T_t\left(-2\be\widehat{Y}^2_s+2|\wh{Y}_s|C\left(|\wh{y}_s|\right)\right)e^{2\be s}ds
+\int^T_t  2|\wh{Y}_s|C |\wh{y}_s|  e^{2\be s}  \ka^{-1} ds +\int^T_t  2|\wh{Y}_s|C |\wh{z}_s|  e^{2\be s}  dL_ {(s-a)^+ }   \nonumber\\
&&-\int^T_t2e^{2\be s}\wh{Y}_s\wh{Z}_sdB_{L_ {(s-a)^+ }   }\nonumber\\
&=&\int^T_t\left(-2\be\widehat{Y}^2_s+2|\wh{Y}_s|C_1\left(|\wh{y}_s|\right)\right)e^{2\be s}ds
+\int^T_t  2|\wh{Y}_s|C |\wh{z}_s|  e^{2\be s}  dL_{(s-a)^+ }    -\int^T_t2e^{2\be s}\wh{Y}_s\wh{Z}_sdB_{L_ {(s-a)^+ }  }\nonumber\\
&\le&\int^T_t\left(-2\be\wh{Y}^2_s+2\frac{C_1}{\sqrt{\la}}|\wh{Y}_s|\sqrt{\la}|\wh{y}_s|
\right)e^{2\be s}ds  + \int_t^T 2\frac{C_1}{\sqrt{\la}}|\wh{Y}_s|\sqrt{\la}|\wh{z}_s| e^{2\be s} d L_{(s-a)^+ }  \nonumber\\
&&-\int^T_t2e^{2\be s}\widehat{Y}_s\widehat{Z}_sd B_{L_{(s-a)^+ } }      \nonumber\\
&\le&\int^T_t\left(-2\be |\wh{Y}|^2_s+2\frac{C_1}{\sqrt{\la}}|\wh{Y}_s|\sqrt{\la}|\wh{y}_s|
\right)e^{2\be s}ds  + \int_t^T 2\frac{C_1}{\sqrt{\la}}|\wh{Y}_s|\sqrt{\la}|\wh{z}_s| e^{2\be s} d L_{(s-a)^+ }
\nonumber\\
&&  -\int^T_t2e^{2\be s}\widehat{Y}_s\widehat{Z}_sd B_{L_{(s-a)^+ } }      \nonumber\\
&\le&\int^T_t\left(\left(-2\be+\frac{ C_1^2}{\la}+\frac {C_1^2}{   \la\ka  }   \right) |\widehat{Y}|^2_s
+\la | \widehat{y} |^2_s \right)e^{2\be s}ds
   + \int_t^T    \la | \widehat{z} |^2_se^{2\be s} dL_{(s-a)^+ }  -\int^T_t2e^{2\be s}\wh{Y}_s\wh{Z}_sd B_{L_{(s-a)^+ } }   \nonumber\\
&=&\int^T_t
 \la  | \widehat{y} |^2_s  e^{2\be s}ds
   + \int_t^T    \la | \widehat{z} |^2_se^{2\be s} dL_{(s-a)^+ }  -\int^T_t2e^{2\be s}\wh{Y}_s\wh{Z}_sd B_{L_{(s-a)^+ } } , 
   \end{eqnarray}
  where
$C_1= (1+\ka^{-1})C$ and  
  $ \la=\frac{C_1^2(\ka +1)  }  {2\be \ka}$. By the same reasoning as that led to \eqref{e:4.9}, we have
$$ 
\E \left[ \sup_{t\in [0, T]} \left| \int_0^t 2e^{2\be s}\wh{Y}_s\wh{Z}_sd B_{L_{(s-a)^+ }} \right| \right]<\infty
$$
and so
$\left\{t\mapsto \int_0^t 2e^{2\be s}\wh{Y}_s\wh{Z}_sd B_{L_{(s-a)^+ } }; \, t\in [0, T]\right\}$
is a martingale.
 
Dropping 1st term on left hand side   of (\ref{esti-norm1}) and taking expectation, we get
\begin{eqnarray}\label{esti-norm2}
 \EE\int^T_0e^{2\be s} |\wh{Z}|^2_sdL_  {(s-a)^+ }  
\leq  \lambda \, \EE   \int^T_0
\bigg(   |\widehat{y}|^2_s  e^{2\be s}ds
   + \int_0^T       | \widehat{z}|^2_se^{2\be s} dL_{(s-a)^+ } \bigg)       
  \le\la  \|(\wh{y},\wh{z})\|^2_{\cM_{a, \be}}.
\end{eqnarray}

It follows from \eqref{esti-norm1} and the Burkholder-Davis-Gundy inequality that  there is a constant $K>0$ so that 
 \begin{eqnarray}\label{esti-norm3}
 &&   \EE \left[ \sup_{t\in [0, T]} \left( e^{2\be t} |\wh{Y}|^2_t  \right) \right]  \nonumber\\
&\le&\la  \|(\wh{y}, \,  \wh{z}) \|^2_{\cM_{a, \be}}
+2\EE \left[ \sup_{t\in [0, T]}\left|\int^T_te^{2\be s}\wh{Y}_s\wh{Z}_sdB_{L_ {(s-a)^+ }  }\right| \right] \nonumber\\
&\le& \la  \|(\wh{y}, \,  \wh{z}) \|^2_{\cM_{a, \be}} 
 +K\EE \left[  \left(\int^T_0e^{4\be s} |\wh{Y}|^2_s  \, |\wh{Z}|^2_sdL_{(s-a)^+ }    \right)^{1/2} \right] \nonumber\\
&\le&  \la  \|(\wh{y}, \,  \wh{z}) \|^2_{\cM_{a, \be}} +K\EE \left[  \left(\sup_{t\in [0, T]} \left(e^{2 \be t}|\wh{Y}_t|^2 \right) \int^T_0e^{2\be s} |\wh{Z}|^2_sd  L_{(s-a)^+ } \right)^{1/2} \right] \nonumber\\
&\le& \la  \|(\wh{y}, \,  \wh{z}) \|^2_{\cM_{a, \be}} +\frac12\EE \left[ \sup_{t\in [0, T]}\left(e^{2\be t} |\wh{Y}|^2_t\right) \right] +\frac{K^2}2 \EE\int^T_0e^{2\be t} |\wh{Z}_s|^2d L_  {(s-a)^+ }   \nonumber\\
&\le& \lambda (1 +K^2) \|(\wh{y},\wh{z})\|^2_{\cM_{a, \be}}+\frac12\EE \left[\sup_{t\in [0, T]}\left(e^{2\be t} |\wh{Y}|^2_t\right) \right],
\end{eqnarray}
where the last inequality follows from (\ref{esti-norm2}). Thus,
\begin{equation}\label{esti-norm4}
\EE \left[ \sup_{t\in [0, T]}\left(e^{2\be t} |\wh{Y}|^2_t\right) \right]
\le   2 \lambda (1 +K^2)   \|(\wh{y},\wh{z})\|^2_{\cM_{a, \be}}.
\end{equation}
Consequently, 
\begin{equation}\label{e:4.11a}
\EE \left[ \int_0^T e^{2\be t} |\wh{Y}|^2_t dt  \right]
\leq T \EE \left[ \sup_{t\in [0, T]}\left(e^{2\be t} |\wh{Y}|^2_t\right) \right]
\le  2   T \lambda (1 +K^2)    \|(\wh{y},\wh{z})\|^2_{\cM_{a, \be}}.
\end{equation}
 Combining (\ref{esti-norm2}) and (\ref{e:4.11a}), we have
\[\|(\wh{Y},\wh{Z})\|^2_{\cM_\be}\le  3  \lambda (T\vee 1)   (1 +K^2)    \|(\wh{y},\wh{z})\|^2_{\cM_{a, \be}} .\]
Taking $\be$ sufficiently large so that 
$$
3  \lambda (T\vee 1)   (1 +K^2)  = \frac{3C_1^2(\ka +1)  }  {2\be \ka}   (T\vee 1)   (1 +K^2)  <1.
 $$
 Then $\Phi$ is a contraction map 
   on ${\mathcal M}[0, T]$ with respect to the Banach norm 
 $\| \cdot \|_{{\mathcal M}_{a, \be} }$. Hence $\Phi$ has a unique fixed point  $(\bar Y, \bar Z)$ in ${\mathcal M}[0, T]$,
 which is the unique solution to  the BSDE  \eqref{y-general} in ${\mathcal M}[0, T]$.   
   \qed

\section{ Stochastic Maximum Principle}\label{S:5}
   In this section,  we establish stochastic maximum principle  using  both spiking and convex variational methods.
 We assume without loss of generality that  $s=0 $ in the state equation\eqref{e:state}.

 \subsection{Spiking variation}

Throughout the remaining of this paper,      we assume the following hypothesis holds. 

\begin{hypothesis}\label{HP1} \rm
  \begin{itemize}
  \item [\rm (1)]
  There exists a constant $L>0$ so that for $\varphi=b,\ \sigma,\ f:  [0,T]\times\RR\times U\to \RR$,
\[|\varphi(t,x,u)-\varphi(t, y, v)|\le L(|x-y|+ |u-v|)
\quad \hbox{and} \quad  |\varphi(t,0,u)|\le L.\]

  \item [\rm (2)]
 $b,\ \sigma $ and $ f$ are  $C^2$ in $x$
 and for $\varphi=b,\ \sigma,\ f $,
$$
|\partial_x\varphi(t,x,u)-\partial_x\varphi(t,  y, v )| +
|\partial^2_x\varphi(t,x,u)-\partial^2_x\varphi(t,  y, v )| 
\le L(|x-y|+ |u-v|) .
$$
\item [\rm (3)]  $h:\ \RR\to\RR$ is  $C^2$ in $x$ and satisfies for any $x, y\in \RR$
 \[|   h_x( x  )-   h_x(  y  )|+ | h_{xx}( x )- h_{xx}( y  )| \le L(|x-y| ) . \]
\end{itemize}
\end{hypothesis}

 In this section, we assume 
without loss of generality that  $s=0 $ in the state equation \eqref{e:state}
  and $a\geq 0$ is the initial value for $R_0$ in \eqref{e:2.1}.
  Suppose that  $\bar u \in \mathcal{U}'_a[0, T] $ is  the optimal control of 
   \begin{equation}\label{e:4.11} 
    J(0, x_0, u^\ast,  a  )  =\inf_{u \in \mathcal{U}'_a[0, T]} J(0,  x_0, u, a  )
 \end{equation}
    and   $\bar x$ is  the corresponding state process.
 Under   Hypothesis \ref{HP1}, by Theorem \ref{T:BSDE},
	 the following  BSDE has a unique solution $(p_t, q_t)$:
   \begin{equation}\label{p}
\left\{\begin{array}{ccl}
d  p(t)&=&-  \left(  b_x (t, \bar x (t), \bar u(t) )p(t) -f_x(t,\bar x (t), \bar u(t)  )   \right)dt-\si_x(t, \bar x(t), \bar u(t))
  q(t) dL_{(t-a)^+ } \\
&&+  q(t) dB_{L_{(t-a)^+}}
\qquad \hbox{for } t\in [0, T], \\
 p (T)&=&-h'( \bar x(T)).
\end{array}\right.
\end{equation}
  Equation \eqref{p} is a backward SDE (see Definition \ref{YZsolution} and Theorem \ref{T:BSDE}),
and its solution is a pair of stochastic processes $\{(p(t), q(t)); t\in [0, T]\}$ 
 that are adapted to the filtration $\{\sF'_t; t\in [0, T]\}$.

 We consider a  second order adjoint process   $(P(t), Q(t))$ determined by  
 \begin{equation}\label{P}
\left\{\begin{array}{ccl}
dP(t)&=&- \left(2  b_x(t, \bar x(t), \bar u(t))P(t) +  b_{xx}(t, \bar x(t), \bar u(t) ) p(t)-f_{xx} (t, \bar x(t), \bar u(t) )     \right)dt \\
&&  - \Big( (  \sigma_x(t, \bar x(t), \bar u(t)))^2P(t)
 +2  \sigma_x(t, \bar x(t), \bar u(t))Q(t),  \\
 && \qquad +\si_{xx} ( t, \bar x(t), \bar u(t)   )q(t)\Big) dL_{(t-a)^+}     +Q(t)dB_{L_{(t-a)^+}} ,\\
P(T)&=&-h''(\bar x(T)).
\end{array}
\right. 
\end{equation}
  By Theorem \ref{T:BSDE},    the above BSDE has a unique solution under {\bf Hypothesis \ref{HP1}}.

To prove the stochastic maximum principle, we need some preparation.
 Firstly, we establish  a moment  estimate, which will be used several times later.

\begin{lemma}\label{ME}
Suppose $Y\in L^2_{\cF}(0,T)$ be the solution to
\begin{equation}\label{LY}
\left\{\begin{array}{ccl}
dY(t)&=&(a(t)Y(t)+\alpha(t))dt+(b(t)Y(t)+\beta(t))dB_{ L_{(t-a)^+}  } ,\\
Y(0)&=&y_0,
\end{array}\right.
\end{equation}
where $|a(t)|,|b(t)|\le A<\infty$ and $a(t)$ and $\beta(t)$ satisfy 
\begin{equation}\label{eq0310d2}
 \int_0^T \left( \left(\EE \left[ |\alpha(s)|^{2k} \right] \right)^{\frac{1}{2k}}
+ \left(\E  \left[|\beta(s)|^{2k} \right] \right)^{\frac{1}{2k}}  \right) ds<\infty  \quad  \hbox{for } k\ge 1.
\end{equation}
Then there is a constant $K>0$ that depends on $k\geq 1$ so that
\begin{eqnarray}\label{eq0310d3}
  \sup_{0\le t\le T}\EE \left[ |Y(t)|^{2k} \right]  &\le& 
K\left( \EE \left[ |Y(0)|^{2k} \right]+ \left(   \int_0^T  \left(\EE \left[ |\alpha(s)|^{2k} \right] \right)^{1/(2k)}ds \right)^{2k} \right.
 \nonumber \\
&& \left. 
+ \left( \left(  \int_0^T   (\EE \left[ |\beta(s)|^{2k} \right] \right)^{1/k}ds \right)^{k}\right).
\end{eqnarray}
\end{lemma}
\pf
First we assume $\alpha(t)$ and $\beta(t)$ are bounded.
For any $ \varepsilon>0$, define $\left<Y(t) \right>_\varepsilon=\sqrt{Y(t)^2+\varepsilon^2}$. By    the It\^{o} formula   formula,
\begin{eqnarray*}
d ( Y^2(t)) &=&2Y(t)dY(t)+d\left<Y\right>_t \\
&=&2Y(t)(a(t)Y(t)+\alpha(t))dt+2Y(t)(b(t)Y(t)+\beta(t))dB_{L_{(t-a)^+} }+(b(t)Y(t)+\beta(t))^2dL{  (t-a)^+  }\\
&=&2Y(t)(a(t)Y(t)+\alpha(t))dt+2Y(t)(b(t)Y(t)+\beta(t))dB_{L_{(t-a)^+} }+(b(t)Y(t)+\beta(t))^2 \kappa^{-1} \1_{\{R_t=0\}} dt  .
\end{eqnarray*}
\begin{eqnarray*}
d\left<Y(t)\right>_\varepsilon&=&
\frac{1}{2}(Y^2(t)+\varepsilon^2)^{-\frac{1}{2}}
(2Y(t)(a(t)Y(t)+\alpha(t))+(b(t)Y(t)+\beta(t))^2   \kappa^{-1} \1_{\{R_t=0\}}       )dt  \\
&&+Y(t)(b(t)Y(t)+\beta(t))(Y^2(t)+\varepsilon^2)^{-\frac{1}{2}}dB_{   L_{(t-a)^+}    } \\
&&-\frac{1}{8}(Y^2(t)+\varepsilon^2)^{-\frac{3}{2}}4Y^2(t)(b(t)Y(t)+\beta(t))^2 \kappa^{-1} \1_{\{R_t=0\}}          dt \\
&=&\frac{1}{2\left<Y(t)\right>_\varepsilon}
(2Y(t)(a(t)Y(t)+\alpha(t))+(b(t)Y(t)+\beta(t))^2    \kappa^{-1} \1_{\{R_t=0\}}\\&&
-\frac{Y^2(t)(b(t)Y(t)+\beta(t))^2 \kappa^{-1} \1_{\{R_t=0\}}   }{\left<Y(t)\right>^2_\varepsilon})dt+\frac{Y(t)(b(t)Y(t)+\beta(t))}{\left<Y(t)\right>_\varepsilon}dB_{L_{ (t-a)^+  }}.
\end{eqnarray*}
\begin{eqnarray*}
d\left<Y(t)\right>^{2k}_\varepsilon&=&
k\left<Y(t)\right>^{2k-2}_\varepsilon
(2Y(t)(a(t)Y(t)+\alpha(t))+(b(t)Y(t)+\beta(t))^2   \kappa^{-1} \1_{\{R_t=0\}}   \\
&&-\frac{Y^2(t)(b(t)Y(t)+\beta(t))^2  \kappa^{-1} \1_{\{R_t=0\}}   }{\left<Y(t)\right>^2_\varepsilon})dt\\
&&+ 2k\left<Y(t)\right>^{2k-2}_\varepsilon
 Y(t)(b(t)Y(t)+\beta(t))  dB_{L_ { (t-a)^+    }}\\
&&+k(2k-1)\left<Y(t)\right>^{2k-4}_\varepsilon Y^2(t)(b(t)Y(t)+\beta(t))^2 \kappa^{-1} \1_{\{R_t=0\}}    dt.
\end{eqnarray*}
Taking expectation,
\begin{eqnarray*}
\EE \left[ \left<Y(t)\right>^{2k}_\varepsilon \right]
&\le&
\EE \left[ \left<Y_0\right>^{2k}_\varepsilon \right]
+C_0\EE\int^t_0\left<Y(s)\right>^{2k-2}_\varepsilon
(Y^2(s)+|\al(s)Y(s)|+\beta^2(s))ds \\
&&+C_0\EE\int^t_0\left<Y(s)\right>^{2k-4}_\varepsilon Y^2(s)(Y^2(s)+\beta^2(s))ds\\
&\le&\EE \left[ \left<Y_0\right>^{2k}_\varepsilon \right] 
+C_1\EE\int^t_0 \left(\left<Y(s)\right>^{2k}_\varepsilon
+|\al(s)|\left<Y(s)\right>^{2k-1}_\varepsilon
+|\beta(s)|^2\left<Y(s)\right>^{2k-2}_\varepsilon \right)ds.
\end{eqnarray*}
 
Let $\varphi_t=\sup_{s\le t} \left(\EE \left[ \left<Y(s)\right>^{2k}_\varepsilon \right] \right)^{\frac{1}{2k}}$
  and $\delta:=\frac{1}{4C_1}$.
 We have by  Young's inequality that for $t\in [0, \delta]$,   
\begin{eqnarray*}
\varphi^{2k}_t&\le&\varphi^{2k}_0+C_1t\varphi^{2k}_t
+C_1\varphi^{2k-1}_t\int^t_0 \left(\EE \left[\alpha^{2k}(s)\right] \right)^{\frac{1}{2k}}ds
+C_1\varphi^{2k-2}_t\int^t_0 \left(\EE \left[ \beta^{2k}(s) \right] \right)^{\frac{1}{k}}ds \\
&\le&\varphi^{2k}_0+\frac{3}{4}\varphi^{2k}_t
+C_2 \left(\int^t_0 \left(\EE \left[ \alpha^{2k}(s) \right] \right)^{\frac{1}{2k}}ds \right)^{2k}
+C_2 \left(\int^t_0 \left(\EE \left[ \beta^{2k}(s) \right] \right)^{\frac{1}{k}}ds \right)^k.
\end{eqnarray*}
Hence 
\[
\varphi^{2k}_t\le 4\varphi^{2k}_0
+C\left( \left(\int^t_0 \left(\EE \left[ \alpha^{2k}(s) \right] \right)^{\frac{1}{2k}}ds \right)^{2k}
+ \left( \int^t_0 \left(\EE \left[ \beta^{2k}(s) \right] \right)^{\frac{1}{k}}ds \right)^k\right).
\] 
Repeating this on $[k\delta, (k+1) \delta]$ for $1\leq k  \leq [T/\delta]$ and $[[T/\delta]\delta,  T]$, where $[a]$ denotes
the largest integer not exceeding $a$, we conclude that there is a  constant $K>0$ that depends on $k\geq 1$ so that
 \begin{eqnarray*}
\sup_{t\in [0, T]} \EE \left[ \left(Y(t)^2+\varepsilon^2\right)^{k}  \right]
&\le& 
K\left( \EE \left[ \left(Y(0)^2+\varepsilon^2\right)^{k}  \right]+ \left(  \int_0^T   \left(\EE \left[ |\alpha(s)|^{2k} \right] \right)^{1/(2k)}ds \right)^{2k} \right.
 \nonumber \\
&& \left. 
+ \left( \left(  \int_0^T   (\EE \left[ |\beta(s)|^{2k} \right] \right)^{1/k}ds \right)^{k}\right).
 \end{eqnarray*}
Taking $\varepsilon\downarrow 0$ yields \eqref{eq0310d3}.   
\qed

\medskip
 
Now, we consider  spiking variation   as in general $\cU'_a[0,T]$ (equivalently, the control domain $U$) may not be convex. 
 Let $\bar u \in \cU'_a [0, T]$ and $\bar x$ be its corresponding state process, that is,  $\bar x$ is the solution to SDE \eqref{e:state}
with $u=\bar u$ and $s=0$. 
  Fix $v(\cdot)\in\cU'_a[0,T]$. Let $E_\ep\in {\cal B} [0,T]$ with Lebesgue measure
$|E_\ep|=\ep$ (for example,  $E_\ep=[\bar{t},\bar{t}+\ep]$). Define
\begin{equation}\label{uep}  
u^\ep(t)=  
\bar u(t) \1_{ E_\ep^c} (t) + 
v(t) \1_{ E_\ep}(t), 
 \end{equation}
 which is an admissible control in $\cU'_a[0,T]$. Let  
\begin{equation}\label{eq0310d5}
\left\{\begin{array}{ccl}
dx^\ep(t)&=&b(t,x^\ep(t),u^\ep(t))dt+\sigma(t,x^\ep(t),u^\ep(t))dB_{L_{  (t-a)^+    }}, \\
x^\ep_0&=&x_0.
\end{array}\right.
\end{equation}

\begin{lemma}\label{pertubx} For any  $T>0$ and integer $k\geq 1$,
$$
\sup_{t\in [0, T]}   \EE \left[ |x^\ep(t)-\bar x(t)|^{2k} \right]   =O(\ep^k).
$$
\end{lemma}

\pf   For $\varphi= b $ or $\sigma$, define
\begin{equation}\label{eq0310d6}
\left\{\begin{array}{ccl}
 \varphi_x(t)&=&\partial_x\varphi(t, \bar x(t), \bar u(t)),\quad
\varphi_{xx}(t) = \partial^2_x\varphi(t, \bar x(t), \bar u(t)), \\
\delta\varphi(t)&=&\varphi(t,\bar x(t),v(t))-\varphi(t, \bar x(t), \bar u(t)), \\
 \delta\varphi_x(t)  &=&   \varphi_x (t, \bar x(t),v(t))-\varphi_x (t, \bar x(t), \bar u(t)),
 \\
 \delta\varphi_{xx}(t) &=&\varphi_{xx} (t, \bar x(t),v(t))-\varphi_{xx} (t, \bar x(t), \bar u(t)).
\end{array}\right.
\end{equation}

 Let $\xi^\ep(t)=x^\ep(t)-\bar x(t)$. Since
\begin{eqnarray*}
&&b(t,x^\ep(t),  u^\ep(t))-b(t,\bar x(t), \bar u(t))\\
&=&(b(t,x^\ep(t), \bar u(t))-b(t,\bar x(t), \bar u(t)))\1_{E_\ep^c}(t)
+(b(t,x^\ep(t),v(t))-b(t, \bar x(t),\bar u(t)))\1_{E_\ep}(t) \\
&=&(b(t,x^\ep(t), \bar u(t))-b(t, \bar x(t), \bar u(t)))\1_{E_\ep^c}(t)
+(b(t,x^\ep(t),v(t))-b(t, \bar x(t), v(t)))\1_{E_\ep}(t)\\
&&+(b(t, \bar x(t),v(t))-b(t, \bar x(t), \bar u(t)))\1_{E_\ep}(t)\\
&=&b(t,x^\ep(t),u^\ep(t))-b(t, \bar x(t),u^\ep(t))
+(b(t,\bar x(t),v(t))-b(t, \bar x(t), \bar u(t)))\1_{E_\ep}(t)\\
&=&b_x(t, \bar x(t)+\theta(x^\ep(t)-\bar x(t)),u^\ep(t))\xi^\ep(t)
+\delta b(t)\1_{E_\ep}(t)\\
&=&\widetilde{b}^\ep_x(t)\xi^\ep(t)+\delta b(t)\1_{E_\ep}(t),
\end{eqnarray*}
we have
\begin{eqnarray*}
d\xi^\ep(t)&=&(b(t,x^\ep(t),u^\ep(t))-b(t, \bar x(t), \bar u(t)))dt
+(\sigma(t,x^\ep(t),u^\ep(t))-\sigma(t, \bar x(t), \bar u(t)))d B_{L_ { (t-a)^+  }  }   \\
&=&(\widetilde{b}^\ep_x(t)\xi^\ep(t)+\delta b(t)\1_{E_\ep}(t))dt
+(\widetilde{\sigma}^\ep_x(t)\xi^\ep(t)+\delta \sigma(t)\1_{E_\ep}(t))dB_{L_ { (t-a)^+  }   }.
\end{eqnarray*}
Integrating over $[0 , t ]$ and taking expectation, by Gronwall's inequality and Burkholder-Davis-Gundy inequality, we have
 \begin{eqnarray}\label{pGBDG}
\sup_{t\in [0, T]}\EE \left[ (\xi^\ep(t))^{2k} \right]
&\le &C_1\{\int_0^t (\EE(\delta b(t)\1_{E_\ep}(t))^{2k})^{\frac{1}{2k}}dt\}^{2k} 
 +C_1\{\int_0^t(\EE(\delta \sigma(t)\1_{E_\ep}(t))^{2k})^{\frac{1}{k}}dt\}^{k}\nonumber\\
&\le &C_2 (\int_0^t\1_{E_\ep}(t)dt)^{2k}
+C_2(\int_0^t\1_{E_\ep}(t)dt)^{k}\nonumber\\
&\le &C_3\ep^{2k}+C_3\ep^{k} 
\le C \ep^{k}.
\end{eqnarray}
This completes the proof.

\medskip

Define $y^\ep(t)$ be the solution to
\begin{eqnarray}\label{yep}\left\{\begin{array}{ccl}
dy^\ep(t)&=&b_x(t)y^\ep(t)dt+(\sigma_x(t)y^\ep(t)+\delta\sigma(t)\1_{E_\ep}(t))dB_{L_{ (t-a)^+  }}  \\
y^\ep_0&=&0.
\end{array}\right. \end{eqnarray}
\begin{lemma}\label{EMy}
  For every $T>0$ and integer $k\geq 1$,
	$\sup_{t\in [0, T]}\EE \left[ |y^\ep(t)|^{2k} \right] =O(\ep^k).$
\end{lemma}
\pf
Integrating (\ref{yep}) over $[0 , t ]$ and taking expectation, by Gronwall's inequality and Burkholder-Davis-Gundy inequality yields the results.
Since the calculation is standard and similar to  (\ref{pGBDG}), we omit it here.

\medskip

\begin{lemma}\label{x-y}
 For every $T>0$ and integer $k\geq 1$,
\[\sup_{t\in [0, T]}\EE \left[ |x^\ep(t)-\bar x(t)-y^\ep(t)|^{2k} \right] =O(\ep^{2k}).\]
\end{lemma}

\pf  Let $\eta^\ep(t)=x^\ep(t)- \bar x(t)-y^\ep(t)=\xi^\ep(t)-y^\ep(t)$. Then
  using Taylor expansion,  
\begin{eqnarray*}
d\eta^\ep(t)&=&\{(b(t,x^\ep(t), \bar u(t))-b(t, \bar x(t), \bar u(t)))\1_{E^c_\ep}(t) \\
&&+(b(t,x^\ep(t),v(t))-b(t,\bar x(t), \bar u(t))) \1_{E_\ep}(t)-b_x(t)y^\ep(t)\}dt \\
&&+\{\sigma(t,x^\ep(t),u^\ep(t))-\sigma(t,\bar x(t), \bar u(t))
-\sigma_x(t)y^\ep(t)-\delta\sigma(t) \1_{E_\ep}(t)\}dB_{L_  { (t-a)^+  }  } \\
&=&\{b_x(t,\bar x(t), \bar u(t))\xi^\ep(t)
+\frac{1}{2}b_{xx}(t,\theta^\ep(t), \bar u(t))(\xi^\ep(t))^2\\
&&+(b(t,x^\ep(t),v(t))-b(t,x^\ep(t), \bar u(t))) \1_{E_\ep}(t)-b_x(t)y^\ep(t)\}dt \\
&&+\{\sigma_x(t)\xi^\ep(t)
+\frac{1}{2}\widetilde{\sigma}^\ep_{xx}(t)(\xi^\ep(t))^2
+\delta\sigma^\ep(t) \1_{E_\ep}(t) \\
&&-\sigma_x(t)y^\ep(t)-\delta\sigma(t) \1_{E_\ep}(t)\}d B_{L_ { (t-a)^+  }   } \\
&=&\{b_x(t)\eta^\ep(t)+\frac{1}{2}\widetilde{b}^\ep_{xx}(t)(\xi^\ep(t))^2
+\delta b^\ep(t) \1_{E_\ep}(t)\}dt   \\
&&+\{\sigma_x(t)\eta^\ep(t)+\frac{1}{2}\widetilde{\sigma}^\ep_{xx}(t)(\xi^\ep(t))^2
+(\delta\sigma^\ep(t)-\delta\sigma(t)) \1_{E_\ep}(t)\}d B_{L_  { (t-a)^+  }   } ,
\end{eqnarray*}
  where $\theta^\ep(t) $ is a state between $x^\ep(t)$  and  $\bar x(t)$, 
  that is,  
  \begin{equation}\label{e:4.26}
  \theta^\ep(t) = \lambda  x^\ep(t) + (1-\lambda) \bar x(t)
  \quad \hbox{for some }   \lambda =\lambda (\omega, t)\in [0, 1].
 \end{equation}
Let $\alpha(t)=\frac{1}{2}\widetilde{b}^\ep_{xx}(t)(\xi^\ep(t))^2
+\delta b^\ep(t) \1_{E_\ep}(t)$ and
$\beta(t)=\frac{1}{2}\widetilde{\sigma}^\ep_{xx}(t)(\xi^\ep(t))^2
+(\delta\sigma^\ep(t)-\delta\sigma(t)) \1_{E_\ep}(t)$.
By Lemma \ref{ME},
\[\sup_{t\in [0, T]}\EE|\eta^\ep(t)|^{2k}
\le C\{\int_0(\EE\alpha^{2k}(s))^{\frac{1}{2k}}ds\}^{2k}
+ C\{\int_0(\EE\beta^{2k}(s))^{\frac{1}{k}}ds\}^{k}.
\]
\begin{eqnarray*}
\EE \left[ \beta^{2k}(t) \right] &\le& C_1\EE|\xi^\ep(t)|^{2k}
+C_1\EE|\delta\sigma^\ep(t)-\delta\sigma(t)|^{2k} \1_{E_\ep}(t)\\
&\le& C_2\ep^{2k}+C_2\EE|x^\ep(t)-\bar x(t)|^{2k} \1_{E_\ep}(t) \\
&\le& C_2\ep^{2k}+C_3\ep^{k} \1_{E_\ep}(t).
\end{eqnarray*}
\[
\left( \int_0^t(\EE \big[ \beta^{2k}(s) \big])^{ {1}/{2k}}ds\right)^k
\le C_4(\int_0^t(\ep^2+\ep \1_{E_\ep}(t))dt)^k
=C_5(\ep^2+\ep^2)^k = C\ep^{2k}.\]
\qed

Define $z^\ep(t)$ be the solution to
\begin{eqnarray}\label{zep}
\left\{\begin{array}{ccl}
dz^\ep(t)&=&(b_x(t)z^\ep(t)+\frac{1}{2}b_{xx}(t)(y^\ep(t))^2
+\delta b(t) \1_{E_\ep}(t))dt \\
&&+(\sigma_x(t)z^\ep(t)+\frac{1}{2}\sigma_{xx}(t)(y^\ep(t))^2
+\delta\sigma_x(t)y^\ep(t) \1_{E_\ep}(t))dB_{L_{ (t-a)^+  }},\\
z^\ep_0&=&0.
\end{array}\right.
\end{eqnarray}

\begin{lemma}
 For every $T>0$ and integer $k\geq 1$,
$\sup_{t\in [0, T]}\EE \left[ |z^\ep(t)|^{2k} \right]=O(\ep^{2k}).$
\end{lemma}

\pf
Integrating (\ref{zep}) over $[0 , t ]$ and taking expectation, by Gronwall's inequality and Burkholder-Davis-Gundy inequality yields the results.
Since the calculation is standard and similar to  (\ref{pGBDG}), we omit it here.

\medskip

\begin{lemma}\label{L:5.7}
  For every $T>0$ and integer $k\geq 1$,
\[\sup_{t\in [0, T]}\EE \left[ |x^\ep(t)-    \bar x(t)  -y^\ep(t)-z^\ep(t)|^{2k} \right] =o(\ep^{2k}). \]
\end{lemma}

\pf   Let     $\xi^\ep(t):=x^\ep(t)-x(t)$,  $\eta^\ep(t):=x^\ep(t)-   \bar x(t)-y^\ep(t) $
and  $\zeta^\ep(t) =\eta^\ep(t)-z^\ep(t)$.   Then  
\[\left\{\begin{array}{ccl}
d\zeta^\ep(t)&=&A(t)dt+D(t)d B_{L_ { (t-a)^+  }  },  \\
\zeta^\ep(0)&=&0.
\end{array}\right.\]
where
\begin{eqnarray*}
A(t)&=&b(t,x^\ep(t),u^\ep(t))-b(t, \bar x(t), \bar u(t))
-b_x(t)y^\ep(t)-b_x(t)z^\ep(t)-\delta b(t) \1_{E_\ep}(t) \\
&&-\frac{1}{2}b_{xx}(t)(y^\ep(t))^2 \\
&=&b(t,x^\ep(t),u^\ep(t))-b(t, \bar x(t), \bar u(t))
-(b(t, \bar x(t),v(t))-b(t, \bar x(t), \bar u(t))) \1_{E_\ep}(t)  \\
&&-b_x(t)y^\ep(t)-b_x(t)z^\ep(t)-\frac{1}{2}b_{xx}(t)(y^\ep(t))^2 \\
&=&b(t,x^\ep(t),u^\ep(t))-b(t, \bar x(t),u^\ep(t))
-b_x(t)(y^\ep(t)+z^\ep(t))-\frac{1}{2}b_{xx}(t)(y^\ep(t))^2 \\
&=&b_x(t, \bar x(t),u^\ep(t))\xi^\ep(t)
+\frac{1}{2}b_{xx}(t,\theta^\ep(t),u^\ep(t))(\xi^\ep(t))^2
-b_x(t)(y^\ep(t)+z^\ep(t)) \\
&&-\frac{1}{2}b_{xx}(t)(y^\ep(t))^2 \\
&=&(b_x(t, \bar x(t),u^\ep(t))-b_x(t, \bar  x(t), \bar u(t)))\xi^\ep(t)
+b_x(t)(\xi^\ep(t)-y^\ep(t)-z^\ep(t))  \\
&&+\frac{1}{2}(b_{xx}(t,\theta^\ep(t),u^\ep(t))
-b_{xx}(t, \bar x(t),u^\ep(t)))(\xi^\ep(t))^2 \\
&&+\frac{1}{2}(b_{xx}(t, \bar x(t),u^\ep(t))
-b_{xx}(t, \bar x(t), \bar u(t)))(\xi^\ep(t))^2
+\frac{1}{2}b_{xx}(t)((\xi^\ep(t))^2-(y^\ep(t))^2) \\
&=&\delta b_x(t) \1_{E_\ep}(t)\xi^\ep(t)+b_x(t)\zeta^\ep(t)
+\frac{1}{2}(b_{xx}(t,\theta^\ep(t),u^\ep(t))
-b_{xx}(t, \bar x(t),u^\ep(t)))(\xi^\ep(t))^2 \\
&&+\frac{1}{2}\delta b_{xx}(t) \1_{E_\ep}(t)(\xi^\ep(t))^2
+\frac{1}{2}b_{xx}(t)((\xi^\ep(t))^2-(y^\ep(t))^2) \\
&=&b_x(t)\zeta^\ep(t)+\alpha^\ep(t),
\end{eqnarray*}
and
\begin{eqnarray*}
D(t)&=&\sigma(t,x^\ep(t),u^\ep(t))-\sigma(t,  \bar x(t), \bar u(t))
-\delta\sigma(t) \1_{E_\ep}(t)
-\delta\sigma_x(t) \1_{E_\ep}(t)y^\ep(t) \\
&&-\sigma_x(t)(y^\ep(t)+z^\ep(t))
-\frac{1}{2}\sigma_{xx}(t)(y^\ep(t))^2 \\
&=&\sigma(t,x^\ep(t),u^\ep(t))-\sigma(t, \bar x(t),u^\ep(t))
-\delta\sigma_x(t) \1_{E_\ep}(t)y^\ep(t)
-\sigma_x(t)(y^\ep(t)+z^\ep(t))\\
&&-\frac{1}{2}\sigma_{xx}(t)(y^\ep(t))^2 \\
&=&\sigma_x(t, \bar x(t),u^\ep(t))\xi^\ep(t)
+\frac{1}{2}\sigma_{xx}(t,\theta^\ep(t),u^\ep(t))(\xi^\ep(t))^2\\
&&-(\sigma_x(t, \bar x(t),u^\ep(t))-\sigma_x(t, \bar  x(t), \bar u(t)))y^\ep(t)
-\sigma_x(t)(y^\ep(t)+z^\ep(t))
-\frac{1}{2}\sigma_{xx}(t)(y^\ep(t))^2 \\
&=&(\sigma_x(t, \bar x(t),u^\ep(t))-\sigma_x(t, \bar x(t), \bar u(t)))(\xi^\ep(t)-y^\ep(t))
+\sigma_x(t)(\xi^\ep(t)-y^\ep(t)-z^\ep(t))  \\
&&+\frac{1}{2}(\sigma_{xx}(t,\theta^\ep(t),u^\ep(t))
-\sigma_{xx}(t, \bar x(t),u^\ep(t)))(\xi^\ep(t))^2 \\
&&+\frac{1}{2}(\sigma_{xx}(t,\bar x(t),u^\ep(t))
-\sigma_{xx}(t, \bar x(t), \bar u(t)))(\xi^\ep(t))^2
+\frac{1}{2}\sigma_{xx}(t)((\xi^\ep(t))^2-(y^\ep(t))^2) \\
&=&\delta\sigma_x(t) \1_{E_\ep}(t)\eta^\ep(t)+\sigma_x(t)\zeta^\ep(t)
+\frac{1}{2}(\sigma_{xx}(t,\theta^\ep(t),u^\ep(t))
-\sigma_{xx}(t,  \bar x(t),u^\ep(t)))(\xi^\ep(t))^2 \\
&&+\frac{1}{2}\delta\sigma_{xx}(t) \1_{E_\ep}(t)(\xi^\ep(t))^2
+\frac{1}{2}\sigma_{xx}(t)((\xi^\ep(t))^2-(y^\ep(t))^2) \\
&=&\sigma_x(t)\zeta^\ep(t)+\beta^\ep(t),
\end{eqnarray*}
with   $\theta^\ep(t)$ as in \eqref{e:4.26} and   
\begin{eqnarray*}
\alpha^\ep(t)&=&\delta b_x(t) \1_{E_\ep}(t)\xi^\ep(t)
+\frac{1}{2}(b_{xx}(t,\theta^\ep(t),u^\ep(t))
-b_{xx}(t,\bar x(t),u^\ep(t)))(\xi^\ep(t))^2 \\
&&+\frac{1}{2}\delta b_{xx}(t) \1_{E_\ep}(t)(\xi^\ep(t))^2
+\frac{1}{2}b_{xx}(t)((\xi^\ep(t))^2-(y^\ep(t))^2),
\end{eqnarray*}
\begin{eqnarray*}
\beta^\ep(t)&=&\delta\sigma_x(t) \1_{E_\ep}(t)\eta^\ep(t)
+\frac{1}{2}(\sigma_{xx}(t,\theta^\ep(t),u^\ep(t))
-\sigma_{xx}(t,\bar x(t),u^\ep(t)))(\xi^\ep(t))^2 \\
&&+\frac{1}{2}\delta\sigma_{xx}(t) \1_{E_\ep}(t)(\xi^\ep(t))^2
+\frac{1}{2}\sigma_{xx}(t)((\xi^\ep(t))^2-(y^\ep(t))^2).
\end{eqnarray*}
By Minkowski inequality and Lemmas  \ref{pertubx}-\ref{x-y}, we have
 \begin{eqnarray*}
  && \int_0^t(\EE \left[  |\alpha^\ep(t)|^{2k} \right])^{\frac{1}{2k}}dt \\
&\le &\int_0^t\Big( \left(\EE \left[ |\delta b_x(t) \1_{E_\ep}(t)\xi^\ep(t)|^{2k} \right] \right)^{\frac{1}{2k}}  
  + \left(\EE \left[ |\frac{1}{2}(b_{xx}(t,    \theta^\ep (t)  ,u^\ep(t))
-b_{xx}(t,   \bar x(t) ,  u^\ep(t)))(\xi^\ep(t))^2|^{2k} \right] \right)^{\frac{1}{2k}} \\
&&+\left(\EE \left[ |\frac{1}{2}\delta b_{xx}(t) \1_{E_\ep}(t)(\xi^\ep(t))^2|^{2k} \right] \right)^{\frac{1}{2k}} 
 + \left(\EE \left[ |\frac{1}{2}b_{xx}(t)((\xi^\ep(t))^2-(y^\ep(t))^2)|^{2k} \right] \right)^{\frac{1}{2k}}\Big) dt \\
 &\le &  C \int_0\Big(  \1_{E_\ep}(t) \left(\EE \left[ |( v(t)-\bar u(t)) \xi^\ep(t)|^{2k} \right] \right)^{\frac{1}{2k}}
+  \left(\EE \left[ |\xi^\ep(t) \, ( \xi^\ep(t))^2|^{2k}\right] \right)^{\frac{1}{2k}}      \\
&&   + \1_{E_\ep}(t) \left(\EE \left[ | (v(t) -\bar u(t))) ( \xi^\ep(t))^2|^{2k}\right] \right)^{\frac{1}{2k}}  
 +\left(\EE \left[ |\frac{1}{2}b_{xx}(t)((\xi^\ep(t))^2-(y^\ep(t))^2)|^{2k} \right] \right)^{\frac{1}{2k}} \Big) dt      \\
  &\le & C \int_0\Big(  \1_{E_\ep}(t)(\EE|\xi^\ep(t)|^{2k})^{\frac{1}{2k}}
  +(    \EE(   |\xi^\ep(t)|     )^{4k})^{\frac{1}{4k}}(\EE|\xi^\ep(t)|^{8k})^{\frac{1}{4k}}+ \1_{E_\ep}(t)(\EE|\xi^\ep(t)|^{4k})^{\frac{1}{2k}}  
 \\
&&+(\EE|\eta^\ep(t)|^{4k})^{\frac{1}{4k}}(\EE|\xi^\ep(t)+y^\ep(t)|^{4k})^{\frac{1}{4k}}\Big) dt \\
&\le&  C \left( \ep^{{3}/{2}}+\ep \cdot \ep^{{1}/{2}} +\ep^2+\ep^{{3}/{2}}\right) \\
&=&o(\ep).
\end{eqnarray*}
In the second inequality above,  we used the  Lipschitz continuity of
 $b_x (t, x, u)$ and $b_{xx}(t, u, x)$ in $x$ for the first  three  terms,
\eqref{e:4.26} for   $\theta^\ep(t)  $  and the definition of  
  $\xi^\ep(t):= x^\ep(t)-\bar x(t) $.  
   The forth inequality holds because of Lemma  \ref{pertubx}  and Lemma  \ref{x-y}.
Similarly, we have
\[\int_0(\EE|\beta^\ep(t)|^{2k})^{\frac{1}{k}}dt=o(\ep^2).\]
  The desired result  then follows   by Lemma  $\ref{ME}$.
\qed

\begin{lemma}\label{Juep}
Assume that  $\bar u$ is an optimal control  and $u^\ep $ is given by (\ref{uep}). Then
\begin{eqnarray*}
J(u^\ep)&=&J( \bar u)+\EE \left[  h'( \bar  x(T))(y^\ep(T)+z^\ep(T)) \right]
+\frac{1}{2}\EE \left[ h''(\bar x(T))(y^\ep(T))^2  \right] \\
&&+\EE\int_0^T \left(  f_x(t)(y^\ep(t)+z^\ep(t))+\frac{1}{2}f_{xx}(t)(y^\ep(t))^2
+\delta f(t) \1_{E_\ep}(t)\right) dt+o(\ep).
\end{eqnarray*}
\end{lemma}

\pf
\begin{eqnarray*}
&&J(u^\ep)-J(\bar u) \\
&=&\EE \left[ h(x^\ep(T))-h(\bar x(T)) \right]
+\EE\int_0^T(f(t,x^\ep(t),u^\ep(t))-f(t,\bar x(t), \bar u(t)))dt\\
&=&\EE \left[  h'( \bar x(T))\xi^\ep(T) \right]
+\frac{1}{2}\EE  \left[ h''(\theta^\ep(T))(\xi^\ep(T))^2 \right] \\
&& 
+\EE\int_0^T \left(f(t,x^\ep(t),u^\ep(t)) -f(t, \bar x(t),u^\ep(t))  +f(t,\bar x(t),u^\ep(t))-f(t, \bar x(t), \bar u(t)) \right)dt\\
&=&\EE \left[  h'(\bar x(T))\xi^\ep(T) \right]
+\frac{1}{2}\EE \left[   h''(\theta^\ep(T))(\xi^\ep(T))^2 \right]
+\EE\int_0^T  \Big( f_x(t, \bar x(t),u^\ep(t))\xi^\ep(t) \\
&&+\frac{1}{2}f_{xx}(t,\theta^\ep(t),u^\ep(t))(\xi^\ep(t))^2\Big) dt
+\EE\int_0 ^T\delta f(t) \1_{E_\ep}(t)dt\\
&=&\EE \left[  h'( \bar x(T))(y^\ep(T)+z^\ep(T)) \right] +\EE  \left[  h'(\bar x(T))\zeta^\ep(T) \right] 
+\frac{1}{2}\EE h''(\bar x(T))(y^\ep(T))^2 \\
&&+\frac{1}{2}\EE \left[ h''(\bar x(T))\eta^\ep(T)(\xi^\ep(T)+y^\ep(T)) \right] 
+\frac{1}{2}\EE \left[  (h''(\theta^\ep(T))-h''(x^\ep(T)))(\xi^\ep(T))^2  \right] \\
&&+\EE\int_0^T  \Big( f_x(t)(y^\ep(t)+z^\ep(t))+f_x(t)\zeta^\ep(t)
+\delta f_x(t) \1_{E_\ep}(t)\xi^\ep(t) \\
&&+\frac{1}{2}(f_{xx}(t,\theta^\ep(t),u^\ep(t))-f_{xx}(t,\bar x(t),u^\ep(t)))(\xi^\ep(t))^2
+\frac{1}{2}\delta f_{xx}(t) \1_{E_\ep}(t)(\xi^\ep(t))^2\\
&&+\frac{1}{2}f_{xx}(t)(y^\ep(t))^2+\frac{1}{2}f_{xx}(t)\eta^\ep(t)(\xi^\ep(t)+y^\ep(t))\Big) dt
+\EE\int_0 ^T\delta f(t) \1_{E_\ep}(t)dt\\
&=&\EE \left[  h'(\bar x(T))(y^\ep(T)+z^\ep(T)) \right]
+\frac{1}{2}\EE \left[  h''(\bar x(T))(y^\ep(T))^2 \right]  \\
&& +\EE\int_0 ^T\Big( f_x(t)(y^\ep(t)+z^\ep(t))    +\frac{1}{2}f_{xx}(t)(y^\ep(t))^2
+\delta f(t) \1_{E_\ep}(t)\Big) dt+o(\ep).
\end{eqnarray*}
\qed

\begin{lemma}\label{L:4.12} 
Let $a\geq 0$, and $p(t)$
 and $y^\ep(t)$  be given by (\ref{p}) and (\ref{yep}), respectively. Then 
  \begin{equation}
\EE \left[ p(T) y^\ep(T)\right] =    \EE\int_0^T \left(f_x(t)y^\ep(t)+  \kappa^{-1} \delta\sigma(t)q(t)  \1_{E_\ep}(t)   \1_{\{R^a_t=0\}} \right) dt  , 
\end{equation}
where $R^a_t :=S_{L_{(t-a)^+}} + a -t$. 
 \end{lemma}

\pf  By It\^o's formula, we have
\begin{eqnarray*}
d(p(t)y^\ep(t))&=&-(b_x(t)p(t)  -f_x(t))y^\ep(t)dt -  \sigma_x(t)q(t)  y^\ep(t) dL_ {(t-a)^+}
+q(t)y^\ep(t)d  B_{L_{(t-a)^+}  }\\
 && +b_x(t)y^\ep(t)p(t)dt+(\sigma_x(t)y^\ep(t)+\delta\sigma(t) \1_{E_\ep}(t))p(t)d   B_{L_ {(t-a)^+}  }\\
&&+q(t) \left(\sigma_x(t)y^\ep(t) dt +\delta\sigma(t) \1_{E_\ep}(t) dL _ {(t-a)^+}   \right)  \\
&=&(f_x(t)y^\ep(t) dt +\delta\sigma(t)q(t) \1_{E_\ep}(t) dL_ {(t-a)^+}   ) \\
&& + \left(q(t)y^\ep(t)+\sigma_x(t)p(t)y^\ep(t) 
 +\delta\sigma(t)p(t) \1_{E_\ep}(t) \right)d  B_{L_ {(t-a)^+}  }   .
\end{eqnarray*}
Thus we have 
 \begin{eqnarray*}
\EE \left[ p(T) y^\ep(T)\right] &=& \EE\int_0^T \left(f_x(t)y^\ep(t) dt +\delta\sigma(t)q(t)  \1_{E_\ep}(t) d L_{(t-a)^+}  \right)  \\
&=&     \EE\int_0^T \left(f_x(t)y^\ep(t)+  \kappa^{-1} \delta\sigma(t)q(t)  \1_{E_\ep}(t)   \1_{\{R^a(t)=0\}} \right) dt  .
\end{eqnarray*}
where we used Proposition \ref{P:3.2} for the the last equality. 
 \qed

\medskip

Recall
\[
\left\{\begin{array}{ccl}
dz^\ep(t)&=&(b_x(t)z^\ep(t)+\frac{1}{2}b_{xx}(t)(y^\ep(t))^2
+\delta b(t) \1_{E_\ep}(t))dt \\
&&+(\sigma_x(t)z^\ep(t)+\frac{1}{2}\sigma_{xx}(t)(y^\ep(t))^2
+\delta\sigma_x(t)y^\ep(t) \1_{E_\ep}(t))d   B_{L_ {(t-a)^+}   }
   \\
z^\ep_0&=&0.
\end{array}\right.
\]
 Using  It\^o's formula and Proposition \ref{P:3.2}, we can derive the following in the same way as  that in the proof of Lemma \ref{L:4.12}.

\begin{lemma}\label{L:4.13} 
Let $p(t)$
 and $z^\ep (t)$ be given by (\ref{p}) and (\ref{zep}) respectively. Then 
  \begin{eqnarray*}
\EE \left[ p(T) z^\ep(T) \right] &=&\EE\int_0^T \Big( f_x(t)z^\ep(t)
+\frac{1}{2}b_{xx}(t)p(t)(y^\ep(t))^2
+  \frac{1}{2\kappa } \sigma_{xx}(t)q(t)(y^\ep(t))^2   \1_{\{R^a(t)=0\}}    \\
&&\hskip 0.5truein + \big(\delta b(t)p(t)+    \kappa^{-1} 
\delta\sigma_x(t)q(t)y^\ep(t)   \1_{\{R^a(t)=0\}  }    \big) \1_{E_\ep}(t)\Big) dt.
\end{eqnarray*}
\end{lemma}

By Lemmas   \ref{Juep}-\ref{L:4.13},    we have
\begin{eqnarray}\label{J-J}
0 &\geq &J(\bar u)-J(u^\ep) \nonumber \\
&=&-\frac{1}{2}\EE \left[ h''( \bar  x(T))(y^\ep(t))^2 \right] 
+\EE\int_0^T  \bigg(\frac{1}{2}\big(  b_{xx} (t) p(t) -f_{xx}(t)   \big)  (y^\ep(t))^2+\big(\delta b(t)p(t) -\delta f(t)  \big)  \1_{E_\ep}(t)\bigg)dt\nonumber\\
&&+\EE\int_0^T  \bigg(\frac{1}{2}\si_{xx} (t) q(t)   (y^\ep(t))^2       + \delta \si (t)q(t)    \1_{E_\ep}(t)\bigg) \kappa^{ -1 } \1_{\{R^a(t)=0\} }  d  t+o(\ep).
\end{eqnarray}

Let $Y^\ep(T)=(y^\ep(T))^2$, then
\begin{eqnarray}\label{Y}
dY^\ep(t)&=&2y^\ep(t)dy^\ep(t)+d\left<y^\ep(t)\right>   \nonumber\\
&=&  2b_x(t) Y^\ep(t) dt
 + \big(   \sigma_x(t)^2   Y^\ep(t) +     (2\sigma_x(t)y^\ep(t)\delta\sigma(t)+(\delta\sigma(t))^2)    \1_{E_\ep}(t)   \big ) d L_{ (t-a)^+   }\nonumber\\
&&+\{2\sigma_x(t)Y^\ep(t)+2y^\ep(t)\delta\sigma(t) \1_{E_\ep}(t)\}dB_{L_{(t-a)^+}}.
\end{eqnarray}

  Recall that $(P(t), Q(t))$ are the solutions for the BSDE \eqref{P}. 
  By It\^o's formula  and Proposition \ref{P:3.2}, we get  the following.

\begin{lemma}\label{PY}
Let $Y^\ep(t)$ and $P(t)$ be given by (\ref{Y}) and (\ref{P}) respectively,  then
 \begin{eqnarray} \label{e:4.25}
 \EE \left[ P(T) Y^\ep(T) \right] &=&\EE\int_0^T\Big( \big((\delta\sigma(t))^2P(t)  \1_{E_\ep}(t) -\si_{xx } (t,x(t),u(t)) q(t)    Y^\ep(t)    \big)    \kappa^{-1}  \1_{\{R^a(t)=0\}}     \nonumber\\
&&\hskip  0.5 truein  -\big (    b_{xx} (t,x(t),u(t)) p(t)-f_{xx} (t,x(t),u(t))            \big)Y^\ep(t) \Big)dt+o(\ep).
\end{eqnarray}
\end{lemma}

  \medskip
   
Fefine a   Hamiltonian
\begin{eqnarray}\label{H1}
  H(t, x, u,  p ):= b(t, x, u) p  -  f(t, x, u)  ,  
\end{eqnarray}
 We have the following stochastic maximum principle.
   
\begin{theorem}[Stochastic maximum principle for spiking variational method]\label{T:4.15}
Suppose that Hypothesis \ref{HP1}  holds.
Let  $(\bar u(\cdot), \bar x(\cdot))$  be a  local optimal pair for the control problem \eqref{e:4.11} with $s=0$
 (in the sense that  for every $v\in \cU_a' [0, T]$, 
$J(\bar u) \leq J(\bar u + \eps v)$ for any $\eps$ with $|\eps|$ sufficiently small). 
  Let     $(p,q)$ and $(P,Q)$ be  the solutions to  \eqref{p} and \eqref{P}, respectively. Then 
 there is a subset ${\cal N}\subset [0, T]$ having zero Lebesgue measure so that for every $t\in [0, T] \setminus {\cal N}$,
     $\bP$-almost surely,  
 \begin{eqnarray}\label{e:5.17}
 &&    H(t, \bar x (t) , \bar u (t),  p (t))   -         H(t, \bar x (t) , u,  p (t))  
  \\
&&  \hskip 0.2truein  - \frac1{\kappa}   \1_{\{R^a(t)=0\}}    \big(\sigma  (t, \bar x (t) , \bar u (t) ) - \sigma  (t, \bar x (t) ,  u \big)  q(t)   
\nonumber \\
&&  \hskip 0.2truein   -\frac1{2\kappa}  \1_{\{R^a(t)=0\}}     ( \sigma  (t, \bar x (t) , \bar u (t) ) - \sigma  (t, \bar x (t) ,  u  )^2 P(t)   
\geq 0  \quad  \hbox{for every } u\in U.        \nonumber   
 \end{eqnarray}
 \end{theorem}

 \pf    Let $v$ be an arbitrary element in $ \cU'_a[0, T]$.
   Recall the definition of $\delta b (t)$ and $\delta \sigma (t)$  from \eqref{eq0310d6}.
Similarly, we define $\delta f (t)= f(t, \bar x(t), v(t)) - f(t,\bar x(t), \bar u(t))$. 
  Substituting \eqref{e:4.25} into (\ref{J-J}) and    applying  Proposition \ref{P:3.2} 
  yields
  \begin{eqnarray}
0 &\geq & J(\bar u)-J(u^\ep) \nonumber \\
& =& \EE\int_0^T    ( \delta b(t ) p(t) -\delta f(t )  ) \1_{E_\ep}(t) dt    +\EE\int_0^T   (\delta \si (t) q(t)+ \tfrac{1}{2}(\delta\sigma(t))^2P(t)  )  \1_{E_\ep}(t)dL_{(t-a)^+}+o(\ep) \nonumber\\
&=&
    \EE\int_0^T      \1_{E_\ep}(t)  \left( ( \delta b(t ) p(t) -\delta f(t )  )   +  \kappa^{-1}
  \left(\delta \si (t) q(t)+ \tfrac{1}{2}(\delta\sigma(t))^2P(t)  \right)  \1_{\{R^a(t)=0\}} \right)  dt  
+o(\ep). \nonumber 
\end{eqnarray}
 Taking $E_\ep=[t,t+\ep]$,  dividing by $\ep$ and then sending $\ep \to 0$,    
 we conclude from the Lebesgue differentiation theorem that for almost every $t\in [0, T]$, 
  \begin{equation}\label{e:4.30}
  \EE   [  \delta b(t ) p(t) -\delta f(t )  ]
  +     \kappa^{-1} \EE  [\1_{\{R^a(t)=0\}}  \left(  \delta \si (t) q(t)+ \tfrac{1}{2}(\delta\sigma(t))^2P(t)   \right)   ]   \leq 0.
\end{equation}
We claim that  for each  $v\in \cU'_a [0, T]$, for almost every $t\in [0, T]$, 
 \begin{equation}\label{e:3.28}
   \delta b(t ) p(t) -\delta f(t )   
  +     \kappa^{-1} \1_{\{R^a(t)=0\}}  \left(  \delta \si (t) q(t)+ \tfrac{1}{2}(\delta\sigma(t))^2P(t)   \right)    \leq 0
  \quad \bP \hbox{-a.s.}
\end{equation}
Suppose the above is not true. 
Then there would be some $v\in  \cU'_a[0, T]$ so that there is a subset $A \subset [0, T]$
having positive Lebesgue measure such that for each $t\in A$,  \eqref{e:3.28} fails on a set of positive $\bP$-measure. 
Denote $\Lambda:= \{(\omega, t)\in \Omega \times [0, T]:  \delta b(t ) p(t) -\delta f(t )   
  +     \kappa^{-1} \1_{\{R^a(t)=0\}}  \left(  \delta \si (t) q(t)+ \tfrac{1}{2}(\delta\sigma(t))^2P(t)   \right) >0 \}$.
  Then $\Lambda$ is an $(\sF'_t)$-progressive measurable set having positive $\bP\times dt$-measure. Define 
  $v^* = v\1_\Lambda + \bar u \1_{\Lambda^c} \in \cU'_a[0, T]$. For this $v^*$, the corresponding 
  $ \delta b(t ) p(t) -\delta f(t )   
  +     \kappa^{-1} \1_{\{R^a(t)=0\}}  \left(  \delta \si (t) q(t)+ \tfrac{1}{2}(\delta\sigma(t))^2P(t)   \right)$
  is non-negative on $\Omega \times [0, T]$ and   is strictly positive on $\Lambda$. 
 This contradiction to the property \eqref{e:4.30} proves the claim \eqref{e:3.28}.
  In particular, for each $u\in U$  \eqref{e:3.28} holds for the deterministic control $v(t)=u$ for  all $t\in [0, T]$.
  Let $U_0$ be a countable dense subset of $U$. Then there is a Borel set ${\cal N}\subset [0, T]$ having zero Lebesgue measure
   so that
   for every $t\in [0, T]\setminus {\cal N}$,  \eqref{e:5.17} holds    for every  $u\in U_0$ almost surely.
  Consequently,   in view of Hypothesis \ref{HP1}, 
  for every  $t\in [0, T]\setminus {\cal N}$,  \eqref{e:5.17} holds    for every  $u\in U $ almost surely.  
    \qed 

\smallskip

\begin{remark} \label{R:4.16} \rm 
   Theorem \ref{T:4.15} 
  can be viewed as a counterpart of the stochastic maximum principle \cite[Theorem 3.3.2]{YZ} for stochastic control driven by sub-diffusions.
Note that when  the L\'evy measure $\nu$ of $S_t$ is zero, that is, when $S_t= \kappa t$, $R^a (t)\equiv  0$ and so $\1_{\{R^a(t)=0\}}=1$.  
  In this case, our result  recovers the classical stochastic maximum principle stated in \cite[Theorem 3.3.2]{YZ}. 
  \end{remark}

\medskip

\subsection{Convex variational method}\label{convex}

In this subsection, we assume  that   the control domain  
 $U\subset \RR $ is convex.   In this case, $\cU_a' [0, T]$ is convex and we are able to derive SMP
 by convex variational method.

\medskip

 Let $\bar u\in \cU_a' [0, T]$.  
  For any $(\varepsilon,v)\in (0,1)\times  \cU'_a [0, T] $, let  
 $   x^{\bar u+\va v }(\cdot) $ be the solutions of
(\ref{e:state}) with $ \bar u+\va v  $. in place of $u$.

\begin{lemma}\label{xpi-x} Suppose Hypothesis  \ref{HP1}  holds. There exits a constant $C$ such that for any $  v\in \cU'_a [0, T]$ and $\va >0$,
  \[    \EE \left[  \sup\limits_{0\leq t \leq T}     |x^{ \bar u+\va v }(t)-x^{\bar u}(t)|^2 \right] \leq C \va ^2     , \]
\end{lemma}

 \pf 
 By Burkholder -Davis-Gundy's maximal inequality,   
\begin{eqnarray*}
 && \EE \left[ \sup\limits_{0\leq t\leq T}   |x^{ \bar u+\va v } (t)  -x^ {\bar u}   (t) |^2  \right] \nonumber\\
& \leq &C_1\EE\Bigg[   \(\int_0^t  |b(s, x^{ \bar u+\va v } (s),  \bar  u(s)+\va v (s)) -b(s,x^ {\bar u} (s), \bar u(s) )| ds\)^2  \nonumber\\
&&+ \( \int_0^t   |\sigma (s, x^{\bar u+\va v } (s), \bar  u(s)+\va v (s)) -\sigma(s,x^{\bar u}(s) ,\bar u(s) )| ^2 d\langle B_{L_{(s-a)^+}}\rangle \)  \Bigg]\nonumber\\
  & \leq &  C_2  \EE \left[ \(  \int_0^t   \left( |   x^{\bar  u+\va v } (s)  -x^ {\bar  u } (s)     |+ \va  | v(s)|   \right)    ds \)^2 \right] \nonumber\\
 && +  C_2   \EE  \left[ \int_0^t \(   |   x^{ \bar  u+\va v } (s)  -x^ {\bar   u }(s)      |+ \va  | v(s)|      \)^2 ds  \right]   \nonumber\\
& \leq &C_3    \EE  \left[   \int_0^t  \sup \limits_{0\leq s\leq t }      |   x^{ \bar  u+\va v } (s)  -x^ {\bar  u } (s)     |^2      ds    +\va ^2   \sup  \limits_{0\leq s\leq t }  |v(s)|^2 \right] .\nonumber \\
\end{eqnarray*}
 The desired inequality now follows from the Gronwall's inequality.  \qed

\medskip
 
Let  $x^{(1)}(t):=x^{1, \bar u, v}(t)$ be the unique solution of

\begin{eqnarray}\label{x1}
\left\{\begin{aligned} d x^{(1)} (t)  =& \,   \big( b_x(t, x^{\bar u}(t), {\bar u}(t)) x^{(1)}(t)   +b_u(t, x^{\bar u}(t), {\bar u}(t)) v(t)\big)dt   \\
 & \, +  \big(\sigma_x(t, x^{\bar u}(t), {\bar u}(t)) x^{(1)}(t)   +\sigma_u(t, x^{\bar u}(t), {\bar u}(t))v (t)\big)    dB_{L_ {(t-a)^+}}, \quad t\in[0,T] ,  \\
x^{(1)}_0  = &  0 .\\
\end{aligned}
\right.
\end{eqnarray}

 \begin{lemma}\label{Esupx1}
 $\displaystyle \EE  \Big[ \sup\limits_{0\leq t\leq T  } | x^{(1)} (t)|^2 \Big]  <\infty. $
  \end{lemma}

 \pf SDE \eqref{x1} can be solved explicitly for $x^{(1)}$.  For simplicity, let  $X_t:=x^{(1)}(t)$, 
 $$ 
 M_t := \int_0^t \sigma_x(r, x^{\bar u}(r), {\bar u}(r)) dB_{L_ {(r-a)^+}} + \int_0^t b_x(r, x^{\bar u}(r), {\bar u}(r)) dr
 $$
 and 
 $$
 Y_t :=  \int_0^t \sigma_u(r, x^{\bar u}(r), {\bar u}(r))v (r)     dB_{L_ {(r-a)^+}} + 
 \int_0^t b_u(r, x^{\bar u}(r), {\bar u}(r)) v(r) dr .
 $$
Then
$$
dX_t = X_t dM_t + dY_t \quad \hbox{with  } X_0=0.
$$
Denote by ${\rm Exp} (- M)$  the Dol\'eans-Dade exponential of the continuous semimartingale $-M$; that is, 
 \begin{equation}\label{e:4.35}
 {\rm Exp} (-M)_t = \exp \left( -M_t -\frac12 \< M\>_t\right).
 \end{equation} 
 Since $d {\rm Exp} (-M)_t= - {\rm Exp} (-M)_t dM_t$, 
   by Ito's formula, 
 \begin{eqnarray*}
 d(  {\rm Exp} (-M)_t X_t) &=&  {\rm Exp} (-M)_t \left( dX_t -X_t dM_t \right) + d \<X,  {\rm Exp} (-M)\>_t \\
 &=& {\rm Exp} (-M)_t dY_t  -  {\rm Exp} (-M)_t X_t d\<M \>_t.
 \end{eqnarray*} 
 Thus 
 $$
 d \left( e^{\<M\>_t } {\rm Exp} (-M)_t X_t \right) =  e^{\<M\>_t  }{\rm Exp} (-M)_t dY_t .
 $$
 It follows that 
 $$
  e^{\<M\>_t } {\rm Exp} (-M)_t X_t  = \int_0^t e^{\<M\>_r }{\rm Exp} (-M)_r dY_r , \quad t\in [0, T].
 $$
 This together with  \eqref{e:4.35}    gives 
 \begin{equation}\label{e:4.36} 
 X_t =    \exp \left( M_t -\frac12 \< M\>_t\right) \int_0^t \exp \left( -M_r +\frac12 \< M\>_r\right) dY_r. 
 \end{equation}
 Note that $\<M\>_t= \int_0^t \sigma_x(r, x^{\bar u}(r), {\bar u}(r)) ^2  dL_{(r-a)^+} \leq \| \sigma_x\|_\infty^2 \kappa^{-1}t$ and 
 for any integer $k\in \mathbbm Z$, 
\begin{eqnarray*}
\E \left[ \exp (kM_t) \right]  
&=& \E  \left[\exp \left( k M_t -   k^2\< M\>_t\right) \exp (k^2\<M\>_t)\right] \\
&\leq&  \left( \E  \left[\exp \left(2 kM_t - 2k^2\< M\>_t\right) \right]\,  \E  \left[\exp (2k^2\<M\>_t)\right] \right)^{1/2}  \\
&\leq & \exp \left(  |k|\| b_x\|_\infty  t +  k^2  \| \sigma_x\|_\infty^2 \kappa^{-1} t\right).
 \end{eqnarray*}
 It then follows from  \eqref{e:4.36}, the Cauchy-Schwarz inequality and the boundedness of  $b_x$, $b_u$, $\sigma_x$ and $\sigma_u$  
 that $ \E \left[ \sup_{t\in [0, T]} |X_t|^2 \right] <\infty$. 
 \qed

\medskip

\begin{lemma}\label{xva2}
 Suppose   {\bf Hypothesis  \ref{HP1}}  holds and let
\[\tilde{x}^{\va}(t)=\frac{x^{\bar u+\va v}(t)-x^{{\bar u} }(t)}{\va }-x^{(1)}(t).\quad    \]
Then
\[
\lim\limits_{\va \rightarrow0} {\EE}  \left[ \sup\limits_{0\leq t\leq T}  |
 \tilde{x}^\va(t)|^2 \right] =0.\]
 \end{lemma}

 \medskip

\noindent{\bf Proof.}     Let
\begin{equation}\label{e:4.37}
x^{1, \va} (t) :=\frac {   x^{\bar u+\va v }(t) -x^{\bar u } (t)       }  {\va}.
\end{equation} 
Then   
\begin{eqnarray}
  d x^{1, \va}(t) 
     &=&\bigg( \int_0^1   b_x\( t, x^{\bar u}(t)+   \la  \( x^ {\bar u+\va v }(t) - {x}^{\bar u}(t) \),  \bar u(t)+
        \va v(t) \)
      d\la \,   x^{1, \va}(t)   \nonumber \\
&& \quad + \int_0^1   b_u\( t,   x^{\bar u}(t), {\bar u}(t)+   \la \va v (t) \) 
    d\la  \,  v(t) \bigg ) dt \nonumber \\
 &&+  \bigg(  \int_0^1   \si_x\( t, x^{\bar u}(t)+   \la  \( x^ {\bar u+\va v }(t) - {x}^{\bar u}(t) \), \bar u(t)+
   \va v(t) 
 \) d\la\,  x^{1, \va}(t)   
   \nonumber \\
&&  \qquad + \int_0^1  \si_u\( t,   x^{\bar u}(t), {\bar u}(t)   +   \la \va v (t) \)     d\la   \, v(t)         \bigg )  dB_{   L_ {(t-a)^+}  } \nonumber \\
  &  =: &   \big( b_x ^\va (t)  x^{1, \va}(t)   +b_u ^\va(t)    v(t)\big)dt    +  \big(\sigma_x^\ep (t)  x^{1, \va}(t)   +\sigma_ u^\va(t)  v (t)\big)    dB_{L_ {(t-a)^+}}.       \label{e:4.38}
\end{eqnarray}
By the same argument as that for Lemma \ref{Esupx1}, we have
\begin{equation} \label{e:4.39}
\sup_{\va \in (0, 1]}  \EE  \Big[ \sup_{t\in [0, T]}   | x^{1, \eps} (t)|^2 \Big]  <\infty
\quad \hbox{and} \quad 
  \EE  \Big[ \sup_{t\in [0, T]}   | \wt x^{\eps} (t)|^2 \Big]  <\infty .
\end{equation}
In view of  \eqref{x1}, 
\begin{eqnarray*} 
d\tilde{x}^\va(t) &=& d \left( x^{1, \va}(t)  -  x^{(1)} (t) \right) \\
&=&   \tilde{x}^\va(t)  \left( b_x ^\va (t)   dt  + \sigma_x^\ep (t)   dB_{L_ {(t-a)^+}} \right)
    \\
   &&+ x^{(1)} (t)  \left( \left(b_x ^\va (t) -  b_x(t, x^{\bar u}(t), {\bar u}(t)) \right) dt   +\left(\sigma_x ^\va (t) -  \sigma_x(t, x^{\bar u}(t), {\bar u}(t)) \right)  dB_{L_ {(t-a)^+}}  \right)  \\
&&   +v(t) \left( \left( b_u ^\va(t)   - b_u(t, x^{\bar u}(t), {\bar u}(t))  \right)  dt    
+  \left( \sigma_ u^\va(t) -\sigma_u (t, x^{\bar u}(t), {\bar u}(t))  \right)    dB_{L_ {(t-a)^+}}   \right)
 \end{eqnarray*}
It follows from    the Burkholder-Davis-Gundy inequality,  Hypothesis  \ref{HP1},  the boundedness of $b_x$, $\sigma_x $,  $b_u$,  $\si _u$,  Lemma \ref{Esupx1} and \eqref{e:4.39} that there are positive constants $C_1$ and $C_2$ so that 
 $$
f(t) := {\mathbb{E}} \left[ \sup_{s\le t}|\tilde{x}^\va(s)|^2 \right]  
 \le C_1\int^t_0f(s)ds+\va^2C_2 \quad \hbox{for every } t\in [0, T]. 
$$
Thus the Gronwall's inequality implies that 
$ f(t)\le \va^2C_2e^{C_1t} $ for every $t\in [0, T]$.
\qed

\medskip

We define the
   the adjoint equation:
\begin{eqnarray}\label{yconv}
\left\{\begin{aligned}  d p(t)   =&  -  \big( b_x(t, x^{\bar u}(t),  {\bar u}(t)) p(t)    -f_x(  t, x^{\bar u}(t),
{\bar u} (t)   ) \big)dt    \\
& \    -  \si_x ( t, x^{\bar u}(t),  {\bar u}(t)  ) q(t)  d L_{(t-a) ^+ }   +q(t)  dB_{L_ {(t-a) ^+ }      }   , \\
p(T)  = &-h _x(x(T)).\\
\end{aligned}
\right.
\end{eqnarray}

\medskip

The following  is the main result of this section.

\begin{theorem}\label{Necessaryconv}
Suppose that  Hypothesis   \ref{HP1}  holds. Let $\bar u(\cdot)\in \mathcal{U}'_a[0, T]  $ 
   be a local  optimal
control  of \eqref{e:4.11} with $a\geq 0$ (in the sense that  for every $v\in \cU_a' [0, T]$, 
$J(\bar u) \leq J(\bar u + \eps v)$ for any $\eps$ with $|\eps|$ sufficiently small) 
and $   \bar x(\cdot)$ be the corresponding  state process.   Then
   for    almost every  $t\in [0,T)$, almost surely
\begin{equation}\label{e:4.31} 
    b_u(t, \bar x (t) , \bar u (t)) p (t) + \kappa^{-1} \sigma_u (t ,\bar x (t), \bar u (t) )  q(t)  \1_{\{R^a_t=0\}} 
 -   f_u(t, \bar x(t) , \bar u(t)  )  =  0.
  \end{equation} 
   \end{theorem}

\noindent{\bf Proof.}
Let $\ep>0$ and $v\in  \cU'_a [0, T]$.
 It follows from Lemma \ref{xva2} that 
\begin{eqnarray}\label{dcostfunc}
 {}\hskip 0.3truein 
 0  & \leq & 
   \lim \limits_{\va  \rightarrow 0}\frac{ J(\bar u(\cdot)+\va v(\cdot) )- J(\bar u(\cdot))     }{\va}
  \\
  &= &  {\mathbb{E}}  [   \int_0^T   (  f_x(t,\bar x (t), \bar u(t))x^{(1)}(t)
+f_u(t, \bar x(t), \bar u(t) ) v(t)  )dt + h_x  (\bar x(T)) x^{(1)}(T)     ].   \nonumber 
\end{eqnarray}
  By Ito's formula,  \eqref{x1} and \eqref{yconv}, 
\begin{eqnarray}\label{x1yito}
  &&  - \EE    [ x^{(1)}(T)  h_x (x(T))  ] 
   \\
 &=&    \EE   \int_0^T   ( x^{(1)}(t)         f_x(  t, \bar x(t),
\bar u(t)   )   +b_u(t, \bar x(t),  \bar u(t))  p(t)v(t)   ) dt   
  \nonumber \\
 && 
+  \EE   \int_0^T  \sigma_u(t, \bar x(t),  \bar u(t)) v (t)       q(t)  dL_ {(t-a) ^+ } .     \nonumber 
\end{eqnarray}
  Substituting  (\ref{x1yito}) into (\ref{dcostfunc}),  we get  
   by Proposition \ref{P:3.2} that 
$$   {\mathbb{E}}     \int_0^T  v(t)  
 \left(f_u(t,\bar x(t), \bar u(t) )- b_u(t,\bar x(t), \bar u(t) ) p(t)          - \sigma_u (t,\bar x(t), \bar u(t) )  q(t)   \kappa^{-1}  \1_{\{R^a_t=0\}} 
 \right) dt  \geq 0  . 
 $$
 Since this holds for any  $v \in \cU'_a [0, T]$, we conclude that 
$$
f_u(t,\bar x(t), \bar u(t) )- b_u(t,\bar x(t), \bar u(t) ) p(t)          - \sigma_u (t,\bar x(t), \bar u(t) )  q(t)   \kappa^{-1}  \1_{\{R^a_t=0\}} =0
$$
$ \bP \times dt $-a.e. on $\Omega \times [0, T]$.   This together with Fubini's theorem establishes  \eqref{e:4.31}.
\qed

 \smallskip
 
\begin{remark} \label{R:4.20}\rm
Using the   Hamiltonian $H(t, x, u, p)$ defined in \eqref{H1}, we can write \eqref{e:4.31} as follows: 
  for almost every $t\in [0, T]$,
$$ 
\frac{\partial}{\partial u} H(t, \bar x(t),  \bar u(t), p(t) )  +     \1_{\{R^a_t=0\}}  \kappa^{-1} \sigma_u (t ,\bar x (t), \bar u (t) )  q(t)  =  0 
\quad \hbox{a.s. }
$$
This is consistent with the SMP \eqref{e:5.17} using spiking variational method,  and can be viewed as its differential version when the control domain is context.
\end{remark}

\section{Sufficient conditions for maximum principle}\label{S:6}  
 
In this section, we assume without loss of generality that  $s=0 $ in the state equation \eqref{e:state},
and denote by $a\geq 0$  the initial value for $R_0$ in \eqref{e:2.1}.

\subsection{General control domain case}

 In view of  Remark  \ref{R:4.20} above and \cite[(3.5.1) and (3.5.2)]{YZ}, 
the following theorem can be regarded as  a  counterpart of \cite[Theorem 3.5.2]{YZ}  for stochastic controls driven by subdiffusion.

\begin{theorem}\label{T:5.1} [Sufficient maximum principle for spiking variational case]
Suppose  Hypothesis \ref{HP1} holds and the function $h(\cdot)$ is convex.
   Fix $x_0 \in \R$ and $a\geq 0$. 
  Let  $\bar u \in  \cU'_a [0, T]$ be an admissible control   and $\bar x:=x^{\bar u, 0, x_0, a} $ be its  state process of \eqref{e:state}.
Let $(p, q)$ be the unique solution for the backward SDE \eqref{p} associated with $(\bar u, \bar x)$.
 Suppose that   for any admissible control $u (t) \in  \cU'_a [0, T]$ and its corresponding state  process $x (t):=x^{u, 0, x_0, a}(t)$,  
   \begin{eqnarray}\label{condi-suffi-1}
&&\EE\int_0^T\bigg(  \big(  b _x(t, \bar x(t),  \bar u(t)) p(t)  -  f_x(t, \bar x(t),  \bar u(t))   \nonumber\\
 && \hskip 0.6truein +  \1_{\{R^a(t)=0\}}   \kappa^{-1} \si_  x(t, \bar x(t),  \bar u(t)) q(t)      \big   ) (x(t)-\bar x(t))    \bigg)dt\nonumber\\
&\geq&\EE \int_0^T  \bigg(   \big(b  (t, x(t), u(t))p(t) -f (t, x(t), u(t))   +
 \1_{\{R^a(t)=0\}}    \kappa^{-1}  \si (t, x(t), u(t)) q(t)   \big) \nonumber \\
&& \quad   - \big ( b  (t, \bar x(t), \bar u(t))p(t) -f (t,\bar x(t), \bar u(t))   + 
 \1_{\{R^a(t)=0\}}   \kappa^{-1} \si (t, \bar x(t), \bar u(t)) q(t)   \big )               \bigg) dt . 
\end{eqnarray}
 Then $(\bar x  (\cdot ),  \bar u ( \cdot)   )$ is an optimal pair 
  for the optimal stochastic control problem \eqref{e:4.11} with $s=0$.  
 \end{theorem}

\pf       Denote $x(t)-\bar x(t)$ by $\eta(t)$. Then $\eta$ satisfies 
   \[\left\{\begin{array}{ccl}
d\eta (t) &=&\big(b_x(t,\bar x(t), \bar u(t)) \eta (t) +\alpha(t)\big) dt+  \big(\si_x(t,\bar x(t),\bar u(t)) \eta (t) +\beta(t)\big)  d B_{L_{(t-a)^+}} , \\
\eta (0) &=&0, 
\end{array}\right.\]
where 
\[\left\{\begin{array}{ccl}
\alpha(t)&=&  - b_x(t,\bar x(t), \bar u(t)) \eta (t) +b(t, x(t), u(t)) -  b (t,\bar x(t), \bar u(t))        \\
\beta(t) &=& - \si_x(t,\bar x(t), \bar u(t)) \eta (t) +\si(t, x(t), u(t)) -  \si (t,\bar x(t), \bar u(t))   .
\end{array}\right.\]
  By Proposition \ref{P:3.2},  
\begin{eqnarray*}
&&\EE \left[  h_x (\bar x(T)  \eta(T)) \right] \\
&=&-\EE \left[ p(T)\eta(T) \right] +\EE \left[ p(0) \eta(0)  \right]   \\
&=&-\EE\int_0^T \big (f_x(t,\bar x (t),  \bar u (t)) \eta (t) +p(t)\al(t) \big) dt  - \EE \int_0^T q(t) \be(t) dL _{ (t-a)^+}    \\
&=&  -\EE\int_0^T  \left( f_x(t,\bar x (t),  \bar u (t)) \eta (t) +p(t)\al(t)  + \kappa^{-1} q(t) \be(t)  \1_{\{R^a(t)=0\}}  \right) dt      \\
&=&\EE\int_0^T\bigg(  \big(  b _x(t, \bar x(t),  \bar u(t)) p(t)  -  f_x(t, \bar x(t),  \bar u(t))  
+     \kappa^{-1}   \si_  x(t, \bar x(t),  \bar u(t)) q(t)    \1_{\{R^a(t)=0\}}      \big   )\eta(t)    \bigg)dt\\
&& -\EE \int_0^T \bigg(  (b(t, x(t), u(t)) -  b (t,\bar x(t), \bar u(t))  ) p(t) + 
    \kappa^{-1}  (\si(t, x(t), u(t))  \nonumber\\
    && \hskip 0.7truein  -    \si (t,\bar x(t), \bar u(t))  ) q(t)      \1_{\{R^a(t)=0\}}    \bigg) dt\\
 &\geq&\EE \int_0^T  \bigg(   \big(b  (t, x(t), u(t))p(t) -f (t, x(t), u(t))   + 
    \kappa^{-1} \si (t, x(t), u(t)) q(t)  \1_{\{R^a(t)=0\}}   \big)  \\
&&  \hskip 0.6truein   - \big ( b  (t, \bar x(t), \bar u(t))p(t) -f (t,\bar x(t), \bar u(t))   +   \kappa^{-1} \si (t, \bar x(t), \bar u(t)) q(t)  \1_{\{R^a(t)=0\}}     \big )               \bigg) dt   \\
&&  -\EE \int_0^T \bigg(  (b(t, x(t), u(t)) -  b (t,\bar x(t), \bar u(t))  ) p(t) + 
  \kappa^{-1}   (\si(t, x(t), u(t))   \nonumber\\
  && \hskip 0.7truein  -  \si (t,\bar x(t), \bar u(t))  ) q(t)   \1_{\{R^a_t=0\}}     \bigg) dt \\ 
&=&-\EE \big(  \int_0^T (  f (t, x(t), u(t)) -f (t,\bar x(t), \bar u(t))      )    \big)    dt,
\end{eqnarray*}
where the inequality is due to  (\ref{condi-suffi-1})
On the other hand, the convexity of $h$ implies 
 \begin{eqnarray*}
    \EE h(x(T)) -\EE h(\bar x(T)) \geq \EE \left[   h_x(\bar x (T) ) ( x(T)-\bar x(T))\right] 
   = \EE \left[  h_x(\bar x (T) ) \xi(T)) \right] . 
 \end{eqnarray*}
 It follows then 
$$
  \EE\left[  \int_0^T (  f (t, \bar x(t), \bar u(t)) dt  + h(\bar x(T)) \right]
  \leq   \EE\left[  \int_0^T (  f (t, x(t), u(t)) dt  + h(x(T)) \right], 
$$
that is, $J(0, x_0, \bar u, a) \leq J(0, x_0, u, a)$ for any admissible control $u\in  \cU'_a [0, T]$
  This proves that $(\bar x  (\cdot ),  \bar u ( \cdot)   )$ is an optimal pair 
  for the  control problem \eqref{e:4.11} with $s=0$.   
\qed

 \begin{remark} Clearly, condition  \eqref{condi-suffi-1} is satisfied if for each $u (t) \in  \cU'_a [0, T]$, 
\eqref{e:3.28} holds for almost every $t\in [0, T]$.
  In particular,  condition\eqref{condi-suffi-1} is satisfied if \eqref{e:5.17} holds.  
\end{remark}

\subsection{Convex control domain case}

In this subsection, we assume the control domain $U$ is convex. In this case, $\cU'_a [0, T]$ is convex.

  In Section \ref{convex}, we have studied first order variation of the state process to get necessary condition.
  Inspired by \cite{lenhartxiongyong2016},    we  now  study its second order variation and use it to derive a sufficient condition for the stochastic maximal principle for the cost functional \eqref{e:2.3} associated with  (\ref{e:state}).
In this subsection, we assume   Hypothesis  \ref{HP1}   holds.

 Suppose $\bar u \in \cU_a' [0, T]$.  For  $(\varepsilon,v')\in (0,1)\times
  \cU_a' [0, T]$, let $v :=v'-\bar u$ and
 $   x^{\bar u+\va v }(\cdot) $ be the solutions of
(\ref{e:state}) with $u$ replaced by $ \bar u+\va v  $.
 
  Let
$ x^{1, \va} (t) :=(  x^{\bar u+\va v }(t) -x^{\bar u } (t)     )/ \va$. 
We know from Lemma \ref{xva2} that $\lim_{\eps \to 0} \E \Big[ \max_{0\leq t\leq T} |x^{1, \va} (t) -x^{(1)}(t) |^2 \Big] =0$,
where the continuous semimartingale $x^{(1)}$ is given by \eqref{x1}.   Define
$ x^{2, \va}(\cdot):= ( x^{1, \va}   (\cdot) -x^{(1)}  (\cdot) )/ \va$.  
 Let $b_x^\va$, $b_u^\va$, $\sigma_x^\va$ and $\sigma_u^\va$ be defined as in \eqref{e:4.38}. 
 We have by  \eqref{x1} and \eqref{e:4.38} that 
\begin{eqnarray}  \label{e:5.2}
 d x^{2, \va} (t)
  &=&\frac 1  {\va }\bigg( \(  b_x ^\va (t)   x^{1, \va}(t)  - b_x(t)x^{(1)}(t)  +  b_u ^\va(t)    v(t) -  b_u(t ) v(t)     \)   \bigg) dt 
  \nonumber  \\
    &&+   \frac 1  {\va }\bigg( \(  \si_x ^\va (t)   x^{1, \va}(t)  - \si_x(t)x^{(1)}(t)  +  \si_u ^\va(t)    v(t) -  \si_u(t ) v(t)     \)   \bigg) d B_{L_{(t-a)^+}} 
    \nonumber  \\
  &=&    \frac 1  \va \(   b_x ^\va (t)   -     b_x   (t)      \)    x^{1, \va}(t)   dt 
  +     \frac 1  \va    b_x (t)    \(  x^{1, \va}(t)  -   x^{(1)}(t)        \)    dt       +  \frac 1 \va  \(   b_u ^\va(t) - b_u(t )       \)   v(t)  dt       \nonumber  \\
   && +  \frac 1  \va \(   \si_x ^\va (t)    -     \si_x   (t)     \)   x^{1, \va}(t)  d B_{L_{(t-a)^+}}   
  +     \frac 1  \va     \si_x (t)  \(    x^{1, \va}(t)      -   x^{(1)}(t)        \)    d B_{L_{(t-a)^+}}     \nonumber  \\
  && +  \frac 1 \va \(  \si_u ^\va(t) - \si_u(t )      \)   v(t)  d B_{L_{(t-a)^+}}  .   
  \end{eqnarray}
    Using the fundamental theorem of calculus, 
\begin{eqnarray*} 
    && \frac 1  \va \(   b_x ^\va (t)   -     b_x   (t)      \)      \\
   &=&  \frac 1 \va  \bigg( \int_0^1   b_x\( t, x^{\bar u}(t)+   \la  \( x^ {\bar u+\va v }(t) - {x}^{\bar u}(t) \), \bar u(t)+\va v(t) \)     d\la  
   - b_x(t,  x^{\bar u}(t), \bar u(t) )         \bigg)    \\
   &=&  \frac 1 \va  \bigg(  \int_0^1   \left( b_x\( t, x^{\bar u}(t)+  
       \la  \,   \eps  x^ {1, \eps}(t) ,    \bar u(t)+\va v(t) \)  
       - b_x\( t, x^{\bar u}(t) , \bar u(t)+\va v(t) \)  \right)  d\la 
    \\
   &&  \qquad +  b_x\( t, x^{\bar u}(t) , \bar u(t)+\va v(t)\)  
        - b_x(t,  x^{\bar u}(t), \bar u(t) )         \bigg)   \\ 
  &=&  \int_0^1  \int_0^1 b_{xx} \(  t, x^{\bar u}(t)+  \theta  \la   \,    \eps  x^ {1, \eps}(t) ,   \bar u(t)+\va v(t)     \) \la d\theta d\la  \, 
   x^ { 1 \ep} (t) \\
  && +    \int_0^1 b_{xv} (t,x^{\bar u} (t) ,    \bar u +\theta \va v ) d\theta v(t)    \\
 &=:&   \,  b_{xx} ^\va(t)  x ^{1, \va}(t)  +  b_{xv} ^{\va}(t)  v(t).
 \end{eqnarray*}
 By  the fundamental theorem of calculus again, 
\begin{eqnarray*} 
      \frac 1 \va  \(   b_u ^\va(t) - b_u(t )       \)   
&=&  \frac 1 \va  \bigg( \int_0^1   b_u\( t,   x^{\bar u}(t), {\bar u}(t)+   \la \va v (t) \)  
    d\la     - b_u(t,  x^{\bar u}(t), \bar u(t) )        \bigg)  \\ 
  &=&  \int_0^1  \int_0^1 b_{uu} \(  t, x^{\bar u}(t),  \bar u(t)+\la \theta \va v(t)     \) \la d\theta d\la  
   v (t) \\
 &=: &    b_{uu} ^\va(t)    v (t)    .   
\end{eqnarray*}
 Similarly, we have
  \begin{eqnarray*} 
     \frac 1  \va \(   \sigma_x ^\va (t)   -     \sigma_x   (t)      \)      
     &=&  \int_0^1  \int_0^1 \sigma_{xx} \(  t, x^{\bar u}(t)+  \theta  \la   \,  \eps  x^ {1, \eps}(t) ,   \bar u(t)+\va v(t)     \) \la d\theta d\la  \, 
   x^ { 1 \ep} (t) \\
  &&+    \int_0^1 \sigma_{xv} (t,x^{\bar u} (t) ,    \bar u +\theta \va v ) d\theta v(t)    \\
 &=:&   \,  \sigma_{xx} ^\va(t)  x ^{1, \va}(t)  +  \sigma_{xv} ^{\va}(t)  v(t),
   \end{eqnarray*}
and
 \begin{eqnarray*} 
      \frac 1 \va  \(   \sigma_u ^\va(t) - \sigma_u(t )       \)   
   =  \int_0^1  \int_0^1 \sigma_{uu} \(  t, x^{\bar u}(t),  \bar u(t)+\la \theta \va v(t)     \) \la d\theta d\la  
   v (t) 
 =:     \sigma_{uu} ^\va(t)    v (t)    .   
\end{eqnarray*}
Thus we have by \eqref{e:5.2},  
\begin{eqnarray*}
 d x^{2, \va} (t)&=&\bigg (  b_{xx} ^\va(t)  (x ^{1, \va}(t))^2  +  b_{xv} ^{\va}(t)  x^{1, \va}(t) v(t)   
 +  b_x(t) x^{2, \va}(t)
 +    b_{uu} ^\va(t)  ( v (t))^2   \bigg)dt\\
 &&+ \bigg ( \si_{xx} ^\va(t)  (x ^{1, \va}(t))^2  +  \si_{xv} ^{\va}(t)  x^{1, \va}(t) v(t)   
 +  \si_x(t) x^{2, \va}(t)
 +   \si_{uu} ^\va(t)  ( v (t))^2   \bigg) dB_{ L_{(t-a)^+}} . 
 \end{eqnarray*}
 By the same argument as that for Lemmas \ref{Esupx1} and \ref{Esupx1}, we have
$$
\sup_{0<\eps \leq 1} \E \left[ \sup_{t\in [0, T]} |x^{2, \va} (t) |^2 \right] <\infty
\quad \hbox{and} \quad  
 \lim_{\va \rightarrow0} {\EE}  \left[ \sup\limits_{  t\in [0, T]}   
   |x^{2, \va} (t) -x^{(2)}(t)  |^2  \right] =0,
   $$
where $x^{(2)}$ is the unique solution for 
\begin{eqnarray}\label{x2-suf} 
d x^{(2)} (t)&=&\bigg ( \frac 1 2 b_{xx}  (t)  (x^{(1)}(t))^2  +  b_{xu}  (t) x^{(1)}(t) v(t)   
 +  b_x(t) x^{(2)}(t)
 +   \frac 1 2 b_{uu}   (t)  ( v (t))^2   \bigg) dt  \\
 &&+ \bigg (\frac 1 2 \si_{xx}   (t)  (x^{(1)}(t))^2  +  \si_{xu}  (t) x^{(1)}(t) v(t)   
 +  \si_x(t) x^{(2)}(t)
 +  \frac 1 2 \si_{uu}  (t)  ( v (t))^2   \bigg) dB_{ L_{(t-a)^+}}   \nonumber 
 \end{eqnarray} 
 with $x^{(2)}(0)=0$.

For $ u \in \cU_a' [0, T]$, we   define the first order and second order variations of the cost functional $J ( u):=J(0, x_0, \bar u, a)$ of \eqref{e:2.3} with $s=0$
as follows. Fix $\bar u\in \cU_a' [0, T]$ and $v\in \cU_a' [0, T]$. For $\eps \in (0, 1]$, define 
\[ J^{1, \va}:= J^{1, \va}(\bar u, v):= \frac{  J(  \bar u + \eps  v  )  - J( \bar u  )     }{ \va   } .\] 
We know from   \eqref{dcostfunc}  that  $J^{(1)}( \bar u,    v) := \lim_{\eps \to 0}  J^{1, \va} (\bar u, v)$ exists and 
\begin{equation}\label{e:5.4} 
J^{(1)}( \bar u,    v)  
={\mathbb{E}} \left[   \int_0^T  \bigg(  f_x(t,\bar x (t), \bar u(t))x^{(1)}(t)
+f_u(t, \bar x(t), \bar u(t) ) v(t) \bigg)dt + h_x  (\bar x(T)) x^{(1)}(T)    \right].
\end{equation}

By a similar argument  as that for $x^{1, \va}$ of \eqref{e:4.37} above,  
 we have  
\begin{eqnarray}\label{J1va} 
  J^{1, \va}( \bar u,     v)     
&=&  \EE \bigg[   \int_0^T \bigg(\int_0^1   f_x\( t, x^{\bar u}(t)+   \la  
   \eps x^{1, \va}(t),    \, 
\bar u(t)+\va v(t) \)   
  d\la \,  x^{1, \va}(t)      \nonumber \\
&&  \qquad +  \int_0^1   f_u\( t,   x^{\bar u}(t), {\bar u}(t)+   \la \va v (t) \) 
      d\la       \,   v(t)               \bigg) dt  \bigg]  \nonumber \\
    && +\EE \left[ \int_0^1  h_x   \( x^{\bar u}(T)+   \la      \eps x^{1, \va}(T)    \)    d\la \,  x^{1, \va}(T)  \right]
    \nonumber \\
 &=: & \,    \EE\left[   \int_0^T  \bigg(  f_x^\va (t ) x^{1, \va}(t)
+f_u^\va (t  ) v(t) \bigg)dt + h_x^\va  ( T ) x^{1, \va} (T)    \right].  
 \end{eqnarray}

Define 
  \[ J^{2, \va}( \bar u,    v) := \frac{  J^{1, \va}( \bar u,   v  )  - J^{(1)}( \bar u  )     }{ \va   }, \] 
 Then by a similar calculation as that for $x^{1, \va}$ and $x^{2, \va}$, we have 
  \begin{eqnarray}\label{J2va}
     J^{2, \va}( \bar u,    v)   
 &= &   \EE\bigg[   \int_0^T  \bigg(  f_{xx}^\va (t ) (x^{1, \va}(t))^2+ 
  f_{xu}^\va (t ) x^{1, \va} (t)  v(t) +f_x (t) x^{2, \va} (t)   
   +f_{uu}^\va (t)( v(t) )^2   \bigg)dt  \nonumber \\
   &&   \qquad + h_{xx}^\va  ( T ) ( x^{1, \va} (T) )^2 
  +h_x(T) x^{2, \va} (T)     \bigg], 
 \end{eqnarray} 
 where
  \begin{eqnarray*}
   f_{xx}^\va (t ) &:=& \int_0^1  \int_0^1 f_{xx}  (  t, x^{\bar u}(t)+  \theta  \la   ( x^ {\bar u+\va v }(t) - {x}^{\bar u}(t) ), \bar u(t)+\va v(t)     ) \la d\theta d\la ,     \\
   f_{xu}^\va (t ) &:=&   \int_0^1 f_{xu}  (  t, x^{\bar u}(t), \bar u(t)+  \theta \va v(t)     )   d\theta  ,      \\
    f_{uu}^\va (t ) &:=& \int_0^1  \int_0^1 f_{uu}  (  t, x^{\bar u}(t) , \bar u(t)+ \theta\la \va v(t)     ) \la d\theta d\la   ,   \\
   h_{xx}^\va (t ) &:=& \int_0^1  \int_0^1 h_{xx}  (  x^{\bar u}(T)+  \theta  \la   ( x^ {\bar u+\va v }(T) - {x}^{\bar u}(T) )    ) \la d\theta d\la .     
   \end{eqnarray*}
Taking  $\va\rightarrow 0$ gives 
\begin{eqnarray}\label{J2}
     J^{(2)  }( \bar u,    v)   &:=&\lim_{\eps \to 0} J^{2, \va} (\bar u, v)  \nonumber \\
 &= &   \EE\bigg[   \int_0^T  \bigg(  f_{xx}(t ) (x^{(1)}(t))^2+ 
  f_{xu} (t )x^{(1)} (t)  v(t) +f_x (t) x^{(2)} (t)+f_{uu}   (t)( v(t) )^2   \bigg)dt  \nonumber \\
  &&  \qquad   + h_{xx}      ( T ) ( x^{   (1)  } (T) )^2 +h_x(T) x^{(2) } (T)       \bigg].  
 \end{eqnarray}  
  It follows immediately that
    $ x^{\bar u+\va v}(t)=x^{\bar u}(t)+\va x^{(1)}(t)+\va^2 x^{(2)}(t)+o(\va^2)$ and 
$J^{\bar u+\va v}(t)=J^{\bar u}(t)+\va J^{(1)}(t)+\va^2 J^{(2)}(t)+o(\va^2)$.
 This establishes the first part of the following theorem.

  \begin{theorem}\label{T:5.2}
 Suppose {\bf Hypothesis \ref{HP1}}  holds. 
   For any $\bar u, v \in  \mathcal{ U }^\prime _a [0,T]$, we have for $\eps \in (0, 1]$, 
\begin{equation}\label{e:5.9}
 J (\bar u+\va v)=J (\bar u)+\va J^{(1)}(\bar u, v)+\va^2 J^{(2)}(\bar u, v)+o(\va^2),
 \end{equation} 
 where  $J^{(1)}(\bar u, v)$ is given by \eqref{e:5.4}
 and  $J^{(2)}(\bar u, v)$ given by \eqref{J2}. Moreover, we have  
  \begin{eqnarray}\label{J1-substi}
     J^{(1)  }( \bar u,    v)     
& = &   \EE     \int_0^T \(  
 f_u(t ) v(t) -b_u(t)p(t)v(t)  \)dt    \nonumber \\
&&     -\EE     \int_0^T     \si_u(t) q(t)v(t)  dL_ {(t-a) ^+ }   , 
  \end{eqnarray}    
and
   \begin{eqnarray}\label{J2-substi}
   J^{(2)  }( \bar u,    v)   
 &= &   \EE  \bigg[  \int_0^T     \bigg(    
 f_{xu} (t )x^{(1)} (t)  v(t) +f_x (t) x^{(2)} (t)+f_{uu}   (t)( v(t) )^2  -\eta(t) 2 x^{(1)}(t) b_u(t)v(t) \nonumber\\
 && \qquad -p(t)b_{xu}  (t) x^{(1)}(t) v(t)    -p(t) b_{uu}(t) (v(t) )^2   -f_x(t) x^{(2)}(t)\bigg)   dt \nonumber\\
 && -\EE  \bigg[  \int_0^T     \bigg( \eta(t) \(  (\si_u(t) )^2 (v(t))^2 +2 \si_x(t)x^{(1)}(t) \si_u(t)v(t)      \)  
  +\gamma(t)2x^{(1)}(t) \si_u(t)v(t) \nonumber\\
   &&   \qquad  +
  q(t) \(\si_{xu}(t)x^{(1)}(t)v(t) +\si_{uu} (t)( v(t) )^2 \) \bigg)  dL_ {(t-a) ^+ }   \bigg] , 
 \end{eqnarray}  
where $x^{(1)}$ is given by (\ref{x1}) ,  $x^{(2)}$ is given by (\ref{x2-suf}) ,
  $(p, q)$ is the unique solution of the BSDE \eqref{p},  
and $\eta$ is   (a part of)   the unique solution of 
the BSDE
\begin{equation}\label{e:5.13}
\left\{\begin{array}{ccl}
d\eta (t) &=&\(f_{xx}(t)-2\eta(t)b_x(t) -p(t) b_{xx}(t) \) dt- \(  \eta(t)\(\si_x(t)\)^2+2 \gamma(t)\si_x(t) -q(t)\si_{xx}(t)     \)  dL_{(t-a)^+}\\&&+ 
  Z_t      B_{L_{(t-a)^+}} , \\
\eta (T) &=& h_{xx}(x(T)).
\end{array}\right.
\end{equation}
   \end{theorem}

 \pf      Identity \eqref{J1-substi} follows directly from \eqref{dcostfunc} and \eqref{x1yito}
 So it remains to show \eqref{J2-substi}.   
  Note that by  Theorem \ref{T:BSDE}, the BSDE \eqref{e:5.13} has a unique solution $(\eta, Z)$ in $\cM [0, T]$.
 
  Using   Ito's formula to $ p(T) x^{(2)}(T) $ and taking expectation yields
\begin{eqnarray*}\label{x2pito}
  &&    \EE   \left[ x^{(2)}(T)  h_x (x(T)) \right] 
= -  \EE   \left[ x^{(2)}(T) p (T) \right]   \nonumber  \\
&=&   -\EE   \int_0^T p  (t)\bigg(     b_{xx}  (t)  (x^{(1)}(t))^2  +  b_{xu}  (t) x^{(1)}(t) v(t)   
 +  b_x(t) x^{(2)}(t)
 +    b_{uu}   (t)  ( v (t))^2            \bigg ) dt \nonumber\\
  &&+\EE        \int_0^T   x^{(2)}(t)   \bigg(    b_x(t) p(t) -f_x(t) \bigg)    dt
  + \EE \int_0^T  x^{(2)}(t)  \sigma_x(t ) q(t)  dL_{(t-a) ^+ }   
 \nonumber \\
&& - \EE  \int_0^T q(t)  \bigg(   \si_{xx}   (t)  (x^{(1)}(t))^2  +  \si_{xu}  (t) x^{(1)}(t) v(t)   
 +  \si_x(t) x^{(2)}(t)
 +   \si_{uu}  (t)  ( v (t))^2        \bigg) dL_ {(t-a) ^+ }  .  \nonumber \\
\end{eqnarray*}
Applying   Ito's formula to $ \eta(t) (x^{(1)}(t))^2 $ and taking expectation gives 
\begin{eqnarray*}\label{x12etaito}
  &&  \EE   \left[( x^{(1)}(T) )^2 h_{xx} (x(T)) \right]  \nonumber \\
&=& -  \EE \int_0^T   \eta (t) 2 x^{(1)}(t) \( b_x(t)x^{(1)}(t) +b_u(t)v(t)  \)     dt   \nonumber \\
  && -    \EE \int_0^T ( x^{(1)}(t) )^2\(f_{xx}(t)-2\eta(t)b_x(t) -p(t) b_{xx}(t)\) dt
 \nonumber  \\
&&   -\EE   \int_0^T  \eta (t) \(     \si_{x}  (t)  x^{(1)}(t)  +  \si_{u}  (t)   v(t)   \) ^2dL_ {(t-a) ^+ }         \nonumber\\
  &&+\EE  \int_0^T  ( x^{(1)}(t)  )^2 \bigg(    \eta (\si_x(t) )^2+2\gamma (t) \si_x(t) +q(t)\si_{xx}(t)  \bigg)    d L_ {(t-a) ^+ }   
 \nonumber \\
&& - \EE  \int_0^T \gamma(t)  2x^{(1)}(t)\bigg(   \si_{x}   (t) (x^{(1)}(t)   +  \si_{u}  (t)   v(t)       \bigg) dL_ {(t-a) ^+ }.  
 \end{eqnarray*}
  Substituting    the above two identities 
    into (\ref{J2}) gives the desired expression \eqref{J2-substi}.
  \qed

\medskip

We end this section with sufficient condition for the stochastic optimal control of \eqref{JSSC-Eq5}. 

\begin{theorem}\label{Sufficiency} 
Suppose $ \bar u \in \cU_a' [0, T]$. 
  If  for every $v\in \cU_a'[0, T]$, $J^{(1)}\( \bar u;v \)= 0$ and 
$ J_2\( \bar u;v \)>0 $. 
Then   $\bar  u(\cdot)$ is a strict local optimal control of \eqref{JSSC-Eq5} 
in the sense that for every $v\in \cU_a' [0, T]\setminus \{0\}$, there is some $\eps_0\in (0, 1]$
so that $J(\bar u) < J(\bar u + \eps)$ for any $0<|\eps | \leq \eps_0$.    
\end{theorem}

\pf   This follows immediately from  \eqref{e:5.9}. 
\qed

\section{Application to a linear quadratic  system}\label{S:7}

  Due to its  simplicity and good structures, there exists many literatures investigating stochastic control problems 
for   linear quadratic (LQ) systems;   see, e.g.,  \cite{wang2015} and the references therein. 
 In this subsection, we apply our main results to a linear quadratic (LQ) system.

Suppose the  state  equation       \eqref{e:state}    is of   the following form:   
\begin{eqnarray*}\label{xE1}
  d x^u (t)  =   ( x^u (t)  +u(t) ) dt +   dB_{L_{(t-a)^+}} \quad \hbox{with } 
x^u (0) =  x_0.  
 \end{eqnarray*}
The objective is to minimize the cost functional 
 \begin{eqnarray*}\label{JE1}
J(u )&:=&  \frac 1 2 \EE \Big[\int_{0}^{T}    u(t)^2 
       dt +   x^u  (T) ^2-2x^u (T) \Big]
\end{eqnarray*}
over $u\in \cU_a' [0, T]$. 
Hence for this model, the control domain $U=\R$ which is convex,  $\sigma (t, x, u)= 1$, $b(t, x, u)=x+u$, $f(t, x, u)= u^2/2$ and $h(x)= (x^2-2x)/2$.
So in  this   case,     the adjoint equation \eqref{p}   
 is of the form 
 \begin{eqnarray}\label{yE1}
\left\{\begin{aligned}
dp(t)  =&  -    p(t)   dt  +q(t)  dB_{L_{(t-a)^+}},  
    \\
p(T)= & 1-x(T) .\\
\end{aligned}
\right.
\end{eqnarray}
By Theorem \ref{Necessaryconv}, 
  if $\bar u \in \cU_a'[0, T]$ is a local optimal control for  $J$ of \eqref{JE1} and $\bar x$ its corresponding state process,   
then 
\begin{eqnarray} \label{uE1}
\bar u(t)=p( t ) \quad \hbox{for } t\in [0, T].
\end{eqnarray}
Without loss of generality, we assume $a<T$.
For  $t\in [a, T]$, we try a solution of (\ref{yE1}) of the form
\begin{eqnarray} \label{solyE1}
p(t)=\phi(t) \bar x  (t)+ \psi(t)   
\end{eqnarray}  
  for some differentiable function $\phi$ and $\psi$ with  $ \phi(T)=-1$  and $ \psi(T)=1$. 
By Ito's formula,  
\begin{eqnarray}
dp(t)&=& \bar x  (t) \phi ' (t) dt+ \phi(t) d  \bar x  (t)+ \psi ' (t) dt \nonumber\\
&=&\bar x  (t) \phi' (t) dt+ \phi(t)    \bar x  (t)  dt+\phi(t)  \bar x   u(t)dt+ \phi(t)   dB_{L_{(t-a)^+} } + \psi ' (t) dt .
 \end{eqnarray}
Comparing the above equation with (\ref{yE1})  and   taking into account of   \eqref{uE1} and  
 (\ref{solyE1}),  we have $\phi (t)= q(t) $ and 
 $$
 - (\phi(t) \bar x  (t)+ \psi(t) ) = \bar x  (t) \phi' (t)  + \phi(t)    \bar x  (t)   +\phi(t)   (\phi(t) \bar x  (t)+ \psi(t) )+    \psi ' (t) 
   $$
   On the time interval $[a, T]$, $\bar x  (t)$ is stochastic.
 Matching the coefficients  of $\bar x  (t)$ as well as  the 0-order term, we get  for $t\in [a, T]$  
\begin{eqnarray}\label{phi}
 \phi(t)'  +2 \phi(t)  +  \phi (t)^2 =0 \end{eqnarray} and
\begin{eqnarray}\label{psi}
 \phi(t)\psi(t)  + \psi(t) ' +\psi(t)=0 .
 \end{eqnarray}
    
  The unique solution to ODE \eqref{phi}   on $[a, T]$  that satisfies   the boundary $  \phi(T)=-1$ is 
\begin{eqnarray}\label{Sphi}
  \phi(t)=  -\frac{2}{e^{2(t-T)}+1}.  
\end{eqnarray}
Putting this into  \eqref{psi} and taking into the account of  the boundary condition $ \psi(T)=1$, we get  for $t\in [a, T]$, 
\begin{eqnarray}\label{Spsi}
\psi(t)=\exp\left(\int_t^T  (\phi(s)+1 )ds\right)
   =    \frac{2 e^{t-T}}{e^{2(t-T)}+1} .   
\end{eqnarray}

On $[0, a]$, $(t-a)^+=0$  and so $B_{L_{(t-a)^+}}=0$. 
  Thus for $t\in [0, a]$, 
  $$
    q(t)=0 , \quad dp(t)  =   -    p(t)   dt    \quad \hbox{and} \quad 
  d \bar x (t) =  \left( \bar x(t)  +u(t)\right) dt =  \left(\bar  x(t)  +p(t)\right) dt.
  $$
  It follows that  for $t\in [0, a]$, $ p(t)=c e^{-t}$ and 
  $$
  \bar x(t)= e^t x_0 + \int_0^t e^{t-s} p(s) ds
  =e^t x_0+   \frac{c }2 ( e^t - e^{-t}). 
  $$
 Since both $\bar x (t)$ and $p(t)$ are continuous processes, they are in particular continuous at $t=a$.
 By    \eqref{uE1}, $p(a)=\phi(a) \bar x  (a)+ \psi(a)$.       It follows that  
 $$
    c e^{-a}= -\frac{2}{e^{2(a-T)}+1}  \left( e^a x_0+   \frac{c }2 ( e^a - e^{-a}) \right)
  +  \frac{2 e^{a-T}}{e^{2(a-T)}+1}.
$$
  Hence 
  $$
  c =  \frac{2( e^{-T}-x_0)}{  e^{-2T} +1}.
  $$
Thus for $t\in [0, a]$, the local optimal control $\bar u (t)$ is given by
\begin{equation}\label{e:6.11}
 \bar u (t) = p(t)=  \frac{2( e^{-T}-x_0)}{  e^{-2T} +1} e^{-t}, 
  \end{equation}
and the corresponding  state process $\bar x (t) $ is given by 
\begin{equation}\label{e:6.12}
\bar x(t)=  e^t x_0+  \frac{   e^{-T}-x_0 }{  e^{-2T} +1} ( e^t - e^{-t}). 
\end{equation}
For $t\in [a, T]$,  plugging \eqref{solyE1} into \eqref{xE1}, we have 
$$
d \bar x (t)  =  \left( 1+ \phi(t) \right) \bar x (t) dt + \psi(t)  dt +   dB_{L_{(t-a)^+}} .
$$
Thus for $t\in [a, T]$, 
$$
  d \left( e^{-\int_a^t \left( 1+ \phi(s) \right) ds }  \bar x(t) \right) =  e^{-\int_a^t \left( 1+ \phi(s) \right) ds } \left( \psi (t) dt + dB_{L_{(t-a)^+}}\right) .
  $$
Consequently, we have by \eqref{Spsi} that  for $t\in [a, T]$, 
\begin{eqnarray}\label{e:6.13}
   \bar x(t) &=&\bar x(a)+     \int_a^t e^{ \int_r^t\left( 1+ \phi(s) \right) ds } \left( \psi (r) dr  + dB_{L_{(r-a)^+}}\right) \nonumber \\
  &=&e^a x_0+  \frac{   e^{-T}-x_0 }{  e^{-2T} +1} ( e^a - e^{-a})   +     \int_a^t \frac{\psi(r)}{\psi (t)} \left( \psi (r) dr  + dB_{L_{(r-a)^+}}\right) \nonumber \\
  &=&    e^a x_0+  \frac{   e^{-T}-x_0 }{  e^{-2T} +1} ( e^a - e^{-a})   
  +    \int_a^t  \frac{e^{r-t} (e^{2(t-T)}+1)}{e^{2(r-T)}+1} 
    \left( \frac{2 e^{r-T}}{e^{2(r-T)}+1}  dr  + dB_{L_{(r-a)^+}}\right) . \nonumber \\
 \end{eqnarray}
     By \eqref{Sphi} and \eqref{Spsi},  the optimal control $\bar u$  on $[a, T]$ is given by 
\begin{equation}\label{e:6.14}
\bar u(t) = p(t)= \phi (t) \bar x(t)+ \psi (t)= \frac{2 e^{t-T} -2\bar x(t) }{e^{2(t-T)}+1}
\quad \hbox{for } t\in [a, T]. 
\end{equation}
With \eqref{e:6.11}-\eqref{e:6.14} and thus explicit expression of the optimal control pair $(\bar u, \bar x)$ in hands, one can compute the value function
$$
V(x_0, a):= \inf_{u\in \cU_a [0, T]} J(0, x_0, u, a)=  J(0, x_0, \bar u, a)= \frac 1 2 \EE\bigg[\int_{0}^{T}    \bar u(t)^2 
       dt + \bar x  (T) ^2-2\bar x(T)\bigg]. 
 $$

  \bigskip
   
 \noindent {\bf Acknowledgement.} The authors thank Jiongmin Yong  and Xunyu Zhou for the helpful discussion on
 the variational inequality (3.3.19)  in \cite{YZ}.  
 They also thank the referees for helpful comments, especially for bringing references 
 \cite{Marcin, NN} to their attention.


\medskip
 
 \small

\begin{thebibliography}{99}

 
\bibitem{Baeumer2009}
 B. Baeumer, M. M. Meerschaert and E. Nane. 
 Brownian subordinators and fractional Cauchy problems. 
 {\it Trans. Amer. Math. Soc.  \bf 361} (2009), 3915-3930.

\bibitem{Ber1}  J. Bertoin.  {\it L\'evy Processes}.  Cambridge University Press, 1996.

\bibitem{Ber2} J. Bertoin.   Subordinators: examples and applications.
{\it Lectures on probability theory and statistics (Saint-Flour, 1997),} 1-91, Lecture Notes in Math., {\bf 1717}.
Springer, Berlin, 1999.

 
\bibitem{Chen2017}
 Z.-Q. Chen.  Time fractional equations and probabilistic representation. 
 {\it Chaos, Solitons and Fractals \bf  102} (2017),  168-174.


 \bibitem{CKKW1} 
 Z.-Q. Chen, P. Kim, T. Kumagai and J. Wang.
 Heat kernel estimates for time fractional equations. 
  {\it Forum Math.  \bf 30} (2018), 1163-1192.   

\bibitem{CKKW2} 
Z.-Q. Chen, P. Kim, T. Kumagai and J. Wang.
 Time fractional Poisson equations: representations and estimates.  
  {\it J. Funct. Anal.   \bf 278}  (2020), no. 2, 108311, 48 pp.   






\bibitem{Oksendal2005}
N. C. Framstad,  B. {\O}ksendal and  A. Sulem.
 A  sufficient stochastic maximum principle for the optimal control of jump diffusions and applications to finance.
 {\it J Optim Theory Appl. \bf 121}  (2004) ,  77-98. 
 
\bibitem{KS} J.  Klafter and I.  M. Sokolov, 
Anomalous diffusion spreads its wings. {\it  Physics World \bf 18} (2005), 29-32.
 

\bibitem{Kushner1965}
H. J. Kushner.
 On the stochastic maximum principle: fixed time of control.
{\it J. Math. Anal. Appl.   \bf 11} (1965), 78-92.

\bibitem{Kushner1972}
H. J. Kushner.
  Necessary conditions for continuous parameter stochastic
optimization problems. {\it SIAM J. Control  Optim. \bf 10} (1972),  550-565.

 
\bibitem{lenhartxiongyong2016}
  S.~Lenhart, J.~Xiong and J.~Yong.
   Optimal controls for stochastic partial differential equations with an application in population modeling.
 {\it  SIAM J. Control Optim. \bf  54} (2016), 495-535.

\bibitem{Li2012} 
J. Li.
 Stochastic maximum principle in the mean-field controls.
 {\it  Automat. \bf 48} (2012),   366-373.

 \bibitem{Marcin}
M. Marcin.  Path properties of subdiffusion-a martingale approach.  {\it Stochastic Models \bf 26} (2),  (2010), 256-271.  


\bibitem{Meerschaert2004}
M. M. Meerschaert and   H.-P.Scheffler.
Limit theorems for continuous-time random walks with infinite mean waiting times.  
  {\it J. Appl.  Probab. \bf  41}   (2004),   623-638 

\bibitem{Meerschaert2014}
 M. M. Meerschaert and P. Straka.
  Semi-Markov approach to continuous time random walk limit processes.  
 {\it Ann. Probab.  \bf  42} (2014), 1699-1723.

\bibitem{MK}  R. Metzler and J. Klafter,  
 The random walk's guide to anomalous diffusion: A fractional dynamics approach.
 {\it  Phys. Rep. \bf 399(1)} (2002), 1-77.


\bibitem{NN}
E. Nane and Y. Ni. A time-changed stochastic control problem and its maximum principle.  
{\it Probab. Math. Statist.   \bf 41}, (2021), 193-215.



\bibitem{Peng1990}
 S. Peng.
  A general stochastic maximum principle for optimal control problems, {\it SIAM J. Control Optim.  \bf  28} (1990),  966-979.


\bibitem{PardouxPeng1990}
E. ~Pardoux and  S. G. Peng.
 Adapted solution of a backward stochastic differential equation.
 {\it Systems \& Control Letters  \bf 14}  (1990), 55-61.

 \bibitem{RY} D. Revuz and M. Yor.
 {\it Continuous Martingales and Brownian Motion}.
 Springer-Verlag, 1991.


\bibitem{Tang1994}
 S. Tang and  X.  Li.
   Necessary conditions for optimal control of stochastic
systems with random jumps. {\it SIAM J. Control Optim.  \bf  32} (1994),
1447-1475.

 \bibitem{wang2015}
G. C. Wang, Z. Wu and J. Xiong.
 A linear-quadratic optimal control problem of forward-backward stochastic differential equations with partial information. 
{\it IEEE Trans. Autom. Control  \bf 60}  (2015), 2904-2916.


  \bibitem{YZ} J. Yong and X. Y. Zhou. 
  {\it Stochastic Controls}. Springer, 1999.

\bibitem{CZ} S. Zhang and Z.-Q. Chen.  
  Fokker-Planck equation for Feynman-Kac transform of anomalous processes.
  {\it  Stoch. Process. Their Appl. \bf 147} (2022), 300-326

\bibitem{CZ2} S. Zhang and  Z.-Q. Chen.  
 Dynamic programming principle and viscosity solutions of  Hamilton-Jacobi-Bellman   equations for sub-diffusions. 
  In preparation.

\bibitem{ZhangChen-FB} S. Zhang and  Z.-Q. Chen. 
Stochastic maximum principle for fully coupled forward-backward stochastic differential equations driven by sub-diffusion. 
Preprint 2023.

\bibitem{zhang2020}
S. Zhang,  X.  Li and J. Xiong.
  Maximum principle for partially observed forward-backward stochastic differential equations with delay.
  {\it System  \& Control Letters \bf 146} (2020),  104812.

 \bibitem{zhang2018}
 S. Zhang, J. Xiong and X. Liu.
   Stochastic maximum principle for partially observed forward-backward stochastic differential equations with jumps and regime switching.   { \it Sci. China Inf. Sci. \bf 61} (2018), 070211:1-070211:13.


 \end{thebibliography}
 \end{document}